\documentclass[12pt,a4paper]{amsart}

\usepackage[margin=1in]{geometry}
\allowdisplaybreaks[1]

\usepackage{amsmath}
\usepackage{amsfonts}
\usepackage{graphicx}
\usepackage{framed}
\usepackage{amssymb}
\usepackage{esvect}
\usepackage{xcolor}
\usepackage{cite}
\usepackage{booktabs}

\usepackage[backref=page]{hyperref}

\renewcommand*{\backref}[1]{}
\renewcommand*{\backrefalt}[4]{%
    \ifcase #1 (Not cited.)%
    \or        (Cited on page~#2.)%
    \else      (Cited on pages~#2.)%
    \fi}

\usepackage{etex}
\usepackage{latexsym}
\usepackage[english]{babel}
\usepackage[utf8x]{inputenc}

\def \R {\mathbb{R}}

\def \e {\varepsilon}
\renewcommand{\epsilon}{\varepsilon}

\usepackage{amsthm}
\theoremstyle{definition}
\newtheorem{definition}{Definition}[section]

\theoremstyle{plain}
\newtheorem{theorem}[definition]{Theorem}

\newtheorem{lemma}[definition]{Lemma}

\renewcommand{\ge}{\geqslant}
\renewcommand{\le}{\leqslant}

\usepackage{mathtools}

\usepackage[initials]{amsrefs}
\numberwithin{equation}{section}

\newcommand{\EQUAZIONE}[1]{\begin{equation}#1\end{equation}}

\begin{document}

\title[Bushfires and Balance: Proactive versus Reactive Policies]{Bushfires and Balance:\\
Proactive versus Reactive Policies\\ in Prescribed Burning}

\begin{abstract}
We introduce a new mathematical model to explore the dynamic relationship between prescribed burning and bushfire occurrence, formulated as a system of ordinary differential equations. The model admits a unique steady-state, and its stability is shown to depend critically on the policy framework governing prescribed burning. In particular, reactive policies,
where prescribed burning is increased in response to bushfire events,
can lead to system instability unless augmented with appropriate feedback control. Conversely, proactive policies that implement moderate, consistent prescribed burning are more effective at stabilizing the system and mitigating the frequency of bushfires. Additionally, higher vegetation regeneration rates contribute positively to equilibrium stability. While this model does not capture the full complexity of fire management, it highlights the potential risks of emotionally driven, reactive responses and underscores the value of preventive, stability-focused strategies.
To the best of our knowledge, this is the first mathematical analysis of preventive burning strategies and their long-term impact on bushfire mitigation.
\end{abstract}

\author[S. Dipierro]{Serena Dipierro}
\address{S. D., 
Department of Mathematics and Statistics,
University of Western Australia,
35~Stirling Highway, WA 6009 Crawley, Australia. }

\email{serena.dipierro@uwa.edu.au}

\author[E. Valdinoci]{Enrico Valdinoci}
\address{E. V., 
Department of Mathematics and Statistics,
University of Western Australia,
35~Stirling Highway, WA 6009 Crawley, Australia. }

\email{enrico.valdinoci@uwa.edu.au}

\author[G. Wheeler]{Glen Wheeler}
\address{G. W.,
School of Mathematics and Applied Statistics,
University of Wollongong,
Northfields Avenue, NSW 2500 Wollongong, Australia. }

\email{glenw@uow.edu.au}

\author[V.-M. Wheeler]{Valentina-Mira Wheeler}
\address{V.-M. W.,
School of Mathematics and Applied Statistics,
University of Wollongong,
Northfields Avenue, NSW 2500 Wollongong, Australia. }

\email{vwheeler@uow.edu.au}

\thanks{Supported by Australian Research Council DP250101080 and DE190100379.}

\subjclass[2020]{34A34; 34D05; 34D20; 34H15; 92-10.}

\keywords{bushfire models; controlled burning; stability analysis; proactive management.}

\maketitle

\section{Introduction}

Prescribed burns are a key strategy for managing bushfire risks by reducing
fuel loads.
The relationship between the frequency of bushfires and the implementation of
prescribed burns is influenced by various factors, such as environmental
conditions and resource availability.
Moreover, as detailed for instance in Section~1.3 of~\cite{REPO}, ``prescribed
burning can be a contentious issue given that management objectives for its
use often conflict with other land management objectives (e.g. effective fuel
management can alter habitat features), and its use can have a range of
off-site consequences (e.g. smoke can cause discomfort and disruption to
affected communities)'', and opposition from communities can rise
``based on a view that prescribed burning is environmentally destructive,
unnecessary, and probably ineffective in reducing bushfire threat''.
Adverse weather conditions, such as extreme heat and dryness, may also limit
the windows for safe prescribed burning, and extreme conditions could lead to
dangerous fire seasons ``regardless of treatment rate'' (see
e.g.~\cite{CLARKE}).

The literature about controlled burning has a very long tradition (see
e.g.~\cite{q-iw0pf294rf8vj30276jgrth}) and prescribed burning was trialled as
early as the 1890s for forest regeneration purposes
(see~\cite{iadjhfw0ireuy03245-63FSp12-35oyu}).
Notably, indigenous peoples worldwide have employed burning practices long
before modern scientific documentation (see e.g.~\cite{BUTZ09,
qa0DOFHWP3REP9H203EPFL9yusdhfiovlf30er-21, amo22}).
In general, the topic of hazard reduction burns is quite controversial and
is often the subject of intense, public debate.
Critics consider it the byproduct of short-term thinking, whose ``benefits are
greater in the immediate aftermath'', but ineffective in the presence of
adverse conditions and potentially dangerous since they may ``activate seed
banks when they burn, which in turn leads to the growth of more flammable
vegetation below'', see~\cite{ojdsacnpreoogub987381i2kwertigf}.
The impact of controlled burning on air quality and biodiversity is also a
subject of environmental concern (see e.g.~\cite{Haikerwal04052015,BIOD}).

Given the complexity of the problem and its practical importance, it is 
desirable to have a quantitative model assessing the relationship between
bushfire intensity and preventive burns.
In terms of available data, determining the variations of prescribed burns
following major bushfire events is challenging.
This is because of variations in reporting and implementation across different
regions and time periods.
As a rule of thumb, however, it seems that recommendations often call for
increased prescribed burning after major bushfires.
This is especially true after major media coverage.
For example, after the 2009 Black Saturday Bushfires the Victorian Bushfires
Royal Commission recommended a doubling of the prescribed burning program,
increasing preventive burning to 5\% of public land annually (see
Recommendation~56 in~\cite{RECOMM}) and this advice was broadly reported in the
media (see e.g.~\cite{m2edia}).
Along the same lines, in the United States, a significant rise in controlled
burns occurred as a response to massive fires in California, Oregon, and other
Western states, such as the catastrophic wildfires in 2017 and 2018 (like the
Camp Fire, which destroyed Paradise, CA).
As a response to these events, in 2020 approximately 9.4 million acres were
treated with prescribed fire in the United States.
The Forest Service aims to treat up to an additional 20 million acres across
the National Forest System over the next decade (see~\cite{2021REP}).\medskip

The scientific literature on wildfires is extensive, particularly with regard to mathematical models and software tools developed to simulate fire front propagation (see e.g.~\cite{BEER, MR4772545, SPARK}
and references therein).

In parallel, a distinct body of research has examined the social dimensions of bushfire prediction and the role of bushfire simulators in policy-making and operational decision-making (see e.g.~\cite{NEA1qeA}).

To the best of our knowledge, however, this work is the first to present a mathematical analysis of preventive burning strategies and to investigate how different long-term management approaches and evolving trends influence the stable mitigation of bushfire events over extended time horizons.
\medskip

The simple model that we present here accounts for the evolution of three
positive quantities: the amount of burning regions in a given area due to
bushfire, the available amount of land which can be potentially set on fire,
and the portion of the area which is burnt intentionally to prevent bushfire.

These quantities are functions of time (though, in a sense, they have to be
thought as ``averaged'' quantities over some interval of time, referring to a given region of reference,
since the
objective here is to understand the effectiveness of a bushfire prevention
plan, which typically ranges over a sufficiently long time elapse).

We will show that the model possesses a unique fixed point.
The stability\footnote{The stability of an equilibrium
in a model of physical relevance is of paramount importance for practical applications: according to~\cite[page 175]{MR2144536},
``one cannot usually pinpoint
positions exactly, but only approximately, so an equilibrium must be stable to
be physically meaningful''.}
of this equilibrium depends on the policy adopted for the
prescribed burning.
Specifically, we will show that when prescribed burning activities are
incremented as a result of bushfire events, the equilibrium is unstable.
Conversely, when the (moderate) prescribed burning replaces natural bushfires, the
equilibrium is stable.

\section{Mathematical setting and main results}

\subsection{Modeling bushfires and prescribed burning}

In our notation, $a$ denotes the amount of unburnt land available at a given
time~$t$, $b$ represents the amount of land subject to preventive controlled
burns, and $f$ accounts for the frequency and intensity of bushfires.

We suppose that the dynamics of this system is regulated by the following
system of three ordinary differential equations.
The first equation encodes the change in available land, which decreases due
to bushfires and controlled burns and increases thanks to a regeneration
factor.
Specifically, we suppose that the available land decreases proportionally with
respect to controlled burns, while the bushfires effect is linear with respect
to both the available land and the bushfire intensity.
As for the regeneration effect, we assume that it is proportional to the
amount of burnt land (i.e., normalising to~$1$ the total amount of land,
to~$1-a$).
These considerations lead to the equation
\EQUAZIONE{\dot a = -\alpha f a-\beta b +\gamma (1-a),}
where~$\alpha$, $\beta$, $\gamma\in(0,+\infty)$.

The second structural equation that we consider describes the change of
bushfire frequency and intensity due to the availability of fuel (assumed to
be proportional to the available land), allowing for a natural decay over
time, namely
\begin{equation}\label{FREQ:EQMD}
\dot f = \zeta a-\eta f,
\end{equation}
where~$\zeta$, $\eta\in(0,+\infty)$.

The term~$-\eta f$ represents the decay of bushfire activity in the absence of available land~$a$, which may result from factors such as fuel depletion or suppression efforts by fire brigades.

The third equation models the variation of controlled burns.
For this, we suppose that there exists a suitable bushfire intensity
threshold~$f_0\in(0,+\infty)$ which represents the ``optimal'' amount
of fires needed in nature for beneficial purposes (such as regenerating the
landscape, removing invasive weeds, releasing nutrients from dead plants and
animals into the soil, creating habitats such as hollows in trees and logs,
and cracking seed coats to trigger germination).

Thus, we suppose that the change in controlled burns is proportional to the
available land, as well as to the excessive bushfire intensity. 
These observations lead to the equation
\begin{equation}\label{PREVEBU}
\dot b = \theta (f-f_0) a,
\end{equation}
where~$\theta\in\R\setminus\{0\}$.

As is customary in the analysis of differential equations
(see e.g.~\cite[page~3]{MR1908418}), our mathematical
framework allows the involved variables to take negative values, in line with the analytic formulation of the system. However, in the context of our application, the quantities of interest,
namely bushfire frequency, available land, and controlled burns,
are assumed to be nonnegative, reflecting their real-world interpretation.

Summarising, the model comprises the system of ordinary differential equations
\begin{equation}\label{ODESY}
\begin{cases}
&\dot a = -\alpha f a-\beta b +\gamma (1-a),\\
&\dot f = \zeta a-\eta f,\\
&\dot b = \theta (f-f_0) a.\end{cases}
\end{equation}
In spite of its simplicity, this system of equations requires a tailored analysis, since
it does not fall into standard categories (for example, no choice of parameters
reduces it to a gradient system).

It is readily seen that the only physically relevant (i.e., corresponding to positive values of the available land~$a$) equilibrium for this system is
\begin{equation}\label{ODE:EQUIL} a_\star:=\frac{\eta f_0}\zeta,
\quad f_\star:=f_0\quad{\mbox{and}}\quad
b_\star:=\frac{\gamma\zeta-\eta f_0(\gamma+\alpha f_0)}{\beta\zeta}
.\end{equation}
This equilibrium is already of practical interest.
The relation between~$a_\star$ and~$f_\star$, that is, ~$f_\star=\frac{\zeta
a_\star}\eta$, quantifies the optimal bushfire intensity threshold in terms of the environmental land.
The third relation in \eqref{ODE:EQUIL}, being equivalent
to~$b_\star=\frac{\gamma}{\beta}-\frac{a_\star (\gamma\eta+\alpha\zeta
a_\star)}{\beta\eta}$, determines the optimal preventive burning given the
environmental land.

For consistency, one expects the available area at equilibrium~$a_\star$ to be
bounded from above by the area of the full region, which has been normalized
to~$1$.
This implies
\begin{equation}\label{LACO1}\eta f_0\le\zeta.\end{equation}
On a similar note, one expects the burning area to be nonnegative and less
than or equal to the full available region, namely~$b_\star\in[0,a_\star]$,
which entails
\begin{equation}\label{LACO2} (\gamma+\alpha f_0)\eta f_0\le\gamma\zeta\le(\beta+\gamma+\alpha f_0)\eta f_0.\end{equation}
When this condition is violated, the prescribed burning strategy becomes physically unfeasible,
resulting either in the complete consumption of the available area through burning or in the absence of any preventive burning activity. In this context, condition~\eqref{LACO2} is particularly informative, as it establishes quantitative thresholds that determine the viability of prescribed burning strategies. Specifically, the parameter~$\zeta$, which reflects the flammability of the available fuel, must lie within a precise range: namely,
it must be sufficiently large to enable effective fuel reduction, yet not so large as to induce uncontrollable burning. These thresholds are explicitly captured by the double inequality in~\eqref{LACO2}, offering practical insight into the constraints under which prescribed burning can be successfully implemented.

The uniqueness of the equilibrium in~\eqref{ODE:EQUIL} identifies the steady
state ensuring safe management of the land.
However, to make this equilibrium ``robust'', it would be desirable that small
deviations do not make the system drift away from the steady state.
That is, we require a stability analysis of the equilibrium, in which the 
sign of the parameter~$\theta$ will play a crucial role -- see Theorems~\ref{THSTA} and~\ref{asdfdvfb20rofjgQLNWedfbVSTA}.\medskip

Regarding the units of measurement in our model, in addition to traditional
land area metrics (e.g., acres), we mention that wildfire frequency can be
expressed in acres per year (total area affected over time), years (average
fire return interval at a location), or fires per unit area per year (a proxy
for regional fire frequency).
Similarly, wildfire intensity can be measured in megawatts, representing the
total power emitted as thermal radiation, often assessed via satellite data (see Table~\ref{TABELL}).
However, we emphasize that the key advantage of our model is its applicability
as a conceptual framework, based on the overall structure of the problem more
than on specific details, providing policymakers with early warnings about
potential structural instabilities in policies set in place without a long-term
plan.

\begin{center}
\begin{table}[h!]\label{TABELL}
\centering
\resizebox{\textwidth}{!}{%
\begin{tabular}{|l|l|l|}
\hline
\textbf{State variables and parameters} & \textbf{Physical dimensions} & \textbf{Methods of measurement} \\
\hline
Land area ($a$, $b$) & Metric units (e.g., acres) & Manual surveys, GPS devices \\
\hline
Wildfire frequency ($f$, $f_0$) & 
\begin{tabular}[c]{@{}l@{}}Acres per year (total area affected over time),\\ years (mean fire return interval), or\\ fires per unit area per year\end{tabular} & Historical records, direct observation \\
\hline
Wildfire intensity (alternative $f$, $f_0$) & Megawatts (thermal radiation power) & Satellite remote sensing \\
\hline
Proportionality parameters & Dimensionally consistent (unit-dependent) & Estimated from empirical data \\
\hline
\end{tabular}
}\caption{\sl State variables and parameters: physical dimensions and methods of measurement}
\end{table}\end{center}

\subsection{Reactive Policy: Increasing prescribed burning after bushfires}
A natural assumption is that the parameter~$\theta$ in~\eqref{PREVEBU}
is positive.
The ansatz in this case is that the more severe the bushfire events, the more
extensive the preventive activity, as an outcome of committee recommendations,
specific policies, increased public sensibility, and media coverage.
Due to this, we consider the case $\theta>0$ corresponding to Reactive Policy.

\begin{figure}[b]
    \centering
    \includegraphics[height=4.59cm]{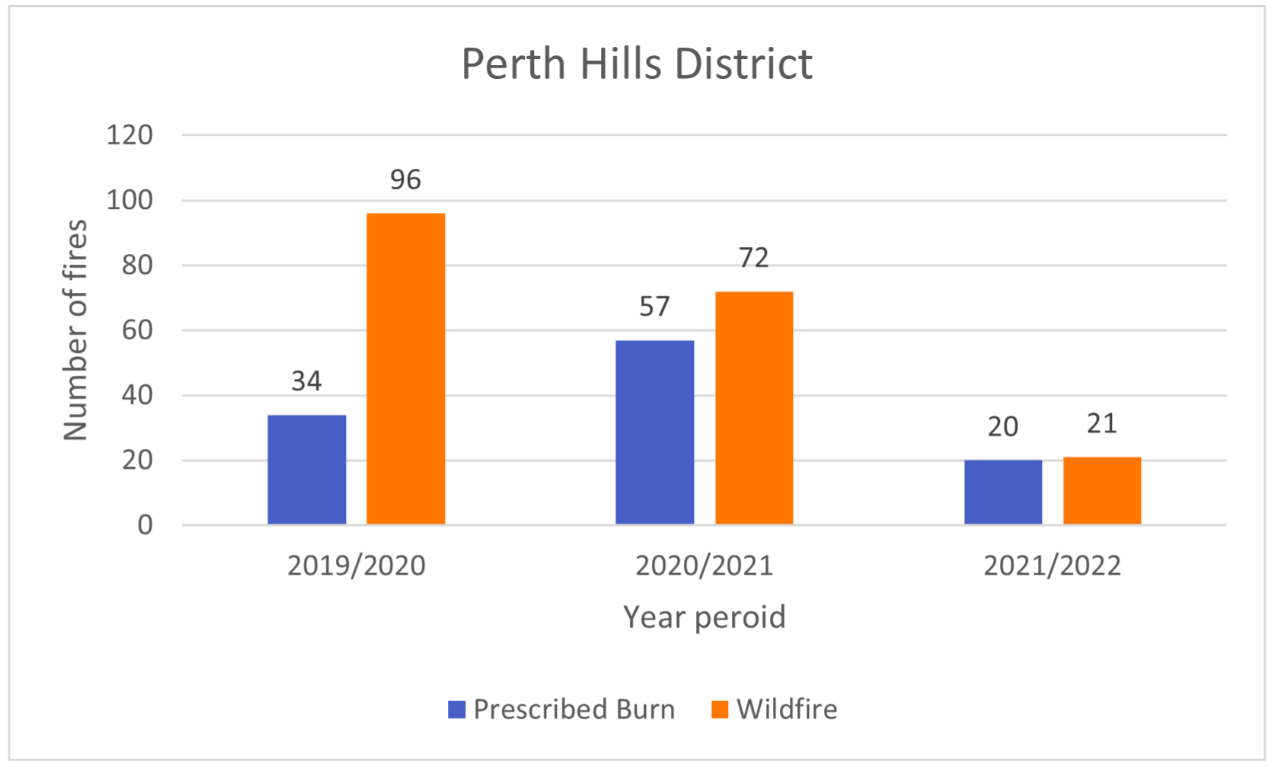}$\,$\includegraphics[height=4.59cm]{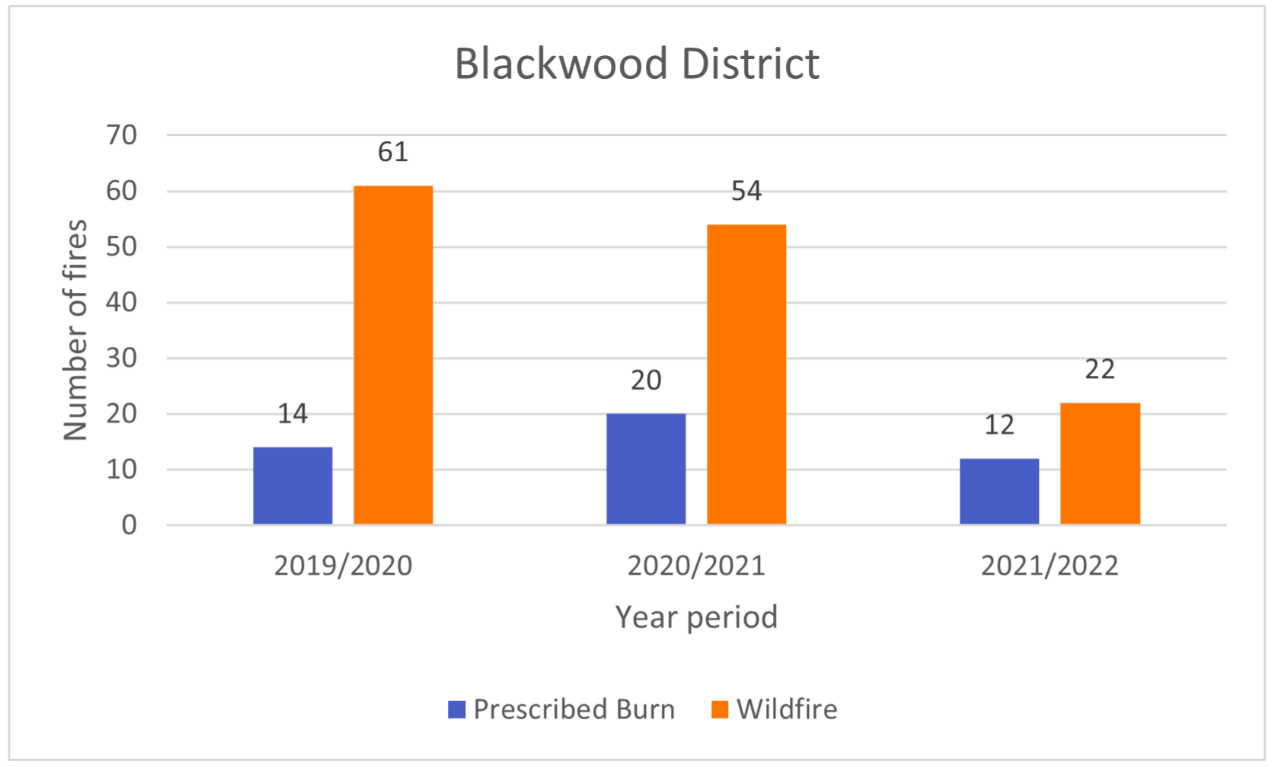}
    \caption{\sl Number of prescribed burns and wildfires in different areas of Western Australia, 2019--2022.
Source:
Brodie James Loller, 2022,
\url{https://storymaps.arcgis.com/stories/6428d63b54e4406eadb62e633f5f41fb}.}
    \label{fig1}
\end{figure}

Interestingly, under this condition, the equilibrium in~\eqref{ODE:EQUIL} is linearly unstable:

\begin{theorem}\label{THSTA}
The equilibrium in~\eqref{ODE:EQUIL} is linearly unstable
when~$\theta\in(0,+\infty)$.

More precisely, depending on the structural parameters~$\alpha$, $\beta$,
$\gamma$, $\zeta$, $\eta$, $\theta$, and~$f_0$,
the Jacobian matrix associated with the
system in~\eqref{ODESY} and evaluated at~\eqref{ODE:EQUIL}
presents: 
\begin{itemize}
\item either three real eigenvalues, of which one is positive and two are negative; or
\item one real, positive eigenvalue and two distinct, complex conjugate
eigenvalues with negative real part.
\end{itemize}
There exist parameters that realize each of the cases above.
\end{theorem}

We observe that, as a byproduct of the eigenvalue classification in
Theorem~\ref{THSTA}, it follows that the unstable space of the equilibrium
in~\eqref{ODE:EQUIL} has dimension one (namely, it is spanned by the only
positive eigenvector).\medskip

For practical purposes, the statement of Theorem~\ref{THSTA}, highlighting the
instability of the only physical equilibrium of the system, could be, in
principle, interpreted as a negative result, hinting at significant concrete
difficulties in managing bushfires merely on the basis of preventive burns.
However, the quantitative analysis and the eigenvalue characterisation in
Theorem~\ref{THSTA} allows one to efficiently utilize an ``augmentation
method'' or ``dynamic extension'': more specifically, one can introduce an
additional ordinary differential equation to the system to stabilize the
equilibrium.
This additional feedback mechanism influences the dynamics and provides
a stabilising damping effect:

\begin{theorem}\label{STAINDO}
There exists a feedback mechanism which, added to~\eqref{ODESY},
makes the equilibrium in~\eqref{ODE:EQUIL} linearly stable.
\end{theorem}

While the statement proposed in Theorem~\ref{STAINDO} may sound vague and
non-specific, the structure provided in Theorem~\ref{THSTA} actually permits an
explicit construction of this feedback mechanism, based on a partial
diagonalisation of the Jacobian matrix and the coupling with a suitable,
linear reaction term:
in this sense, the proof of Theorem~\ref{STAINDO} is more interesting than its
statement, since it exhibits a ``constructive'' procedure to stabilize the
fixed point.
Indeed, from a mathematical perspective, this feedback mechanism penalizes the
positive eigenvalue direction, making its identification useful for wildfire
management.
Furthermore, recognising the potential risk associated with instability can
help justify additional preventive measures to regulate fuel availability.

\begin{figure}[b]
    \centering
    \includegraphics[height=10.59cm]{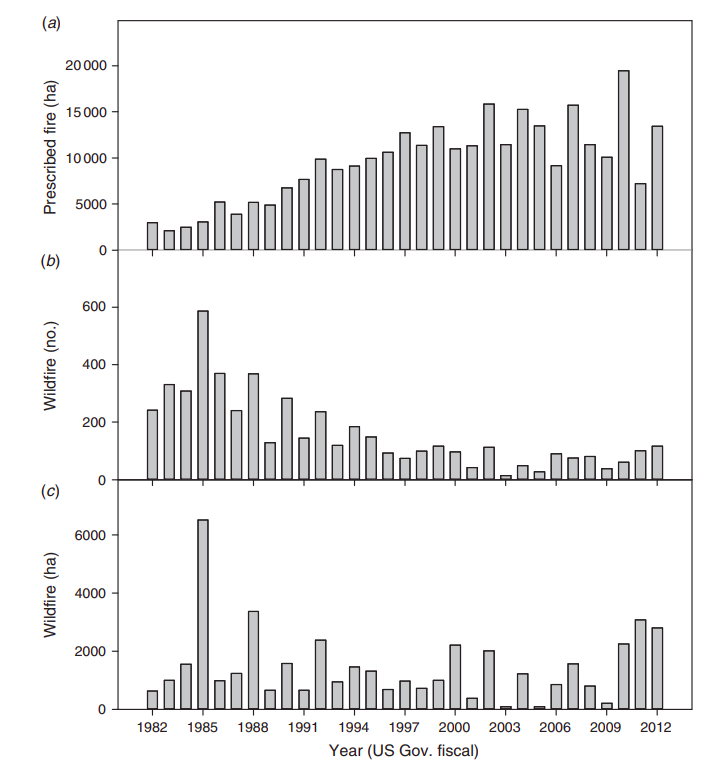}
    \caption{\sl Trends in (a) prescribed fire area (hectares burned), (b) wildfire incidence (number of
unplanned ignitions) and (c) wildfire area (hectares burned by wildfire) from 1982 to 2012 at Fort
Benning, GA.
Source: Figure~1 in~\cite{iwfkhvgri0eotAsd}.}
    \label{fig3}
\end{figure}

\subsection{Proactive Policy: Increasing (with moderation) prescribed burning in the absence of
bushfires}
Interestingly, another method to stabilize the equilibrium in~\eqref{ODE:EQUIL}
consists in modifying the strategy of the preventive burning.
For example, one could perform less controlled burning in the presence of a
highly intense natural fire activity.
The ansatz in this case is that the severe bushfire events have already
removed the available fuel from the soil, hence less preventive activity is
necessary.
Low bushfire activity leaves excessive fuel in the land which needs to be preventively removed via controlled burning.

Bushfires may be the outcome of extreme weather conditions, making controlled
burns unsafe or impossible, which may also suggest that a proactive policy
needs to decrease prescribed burning when bushfires are active.

This controlled burning strategy would maintain the form of~\eqref{PREVEBU},
changing only the sign of its right-hand-side, leading to
\EQUAZIONE{\dot b = -\vartheta (f-f_0) a,}
with~$\vartheta:=-\theta\in(0,+\infty)$.

{F}rom the mathematical standpoint, this plan leaves the system of equations in~\eqref{ODESY}
formally unchanged, except that now~$\theta=-\vartheta$ is negative
(we stress that the equilibrium in~\eqref{ODE:EQUIL}
and conditions~\eqref{LACO1} and~\eqref{LACO2}
do not involve the parameter~$\theta$ and hence remain unaffected by its sign).

\begin{figure}[t]
    \centering
    \includegraphics[height=7.59cm]{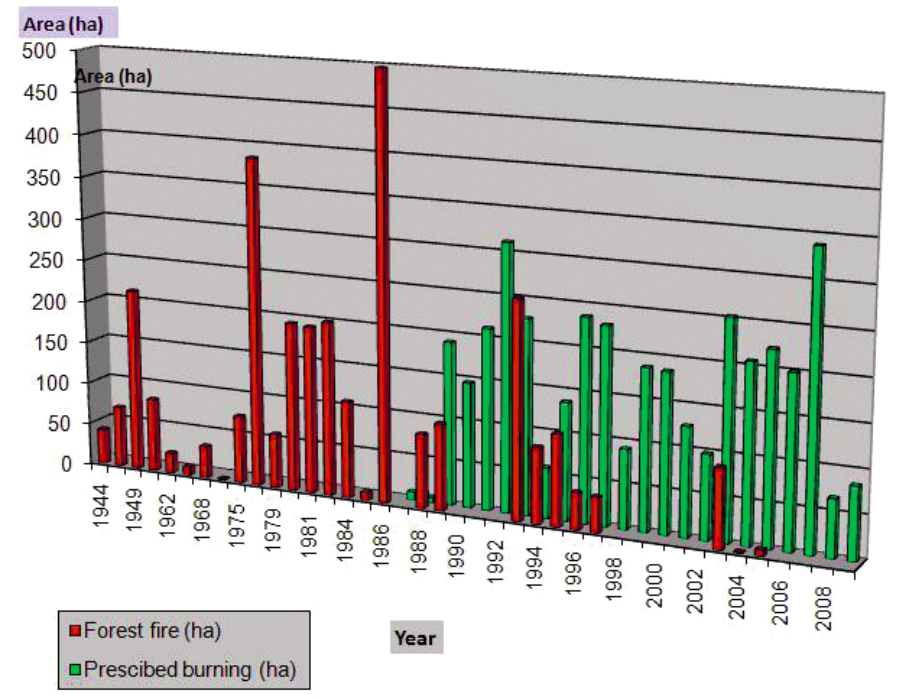}
    \caption{\sl The evolution of forest fires and prescribed burning in the Massif des Madres.
Source: Figure~9 in~\cite{LAMBE}.}
    \label{fig2}
\end{figure}

The relevant fact is that, with this program of preventive fires,
the equilibrium \eqref{ODE:EQUIL} becomes stable for
suitable parameter ranges:

\begin{theorem}\label{asdfdvfb20rofjgQLNWedfbVSTA}
The equilibrium in~\eqref{ODE:EQUIL} is linearly stable for the system in~\eqref{ODESY}
when~$\theta=-\vartheta\in(-\infty,0)$, as long as
\begin{equation}\label{prim:ca}
\left( \alpha f_0 + \gamma + \eta \right) \left( 2 \alpha f_0 \eta + \gamma \eta \right)> \beta \eta \vartheta f_0.
\end{equation}
  
More precisely, under assumption~\eqref{prim:ca},
depending on the structural parameters~$\alpha$, $\beta$,
$\gamma$, $\zeta$, $\eta$, $\theta$, and~$f_0$, the
Jacobian matrix associated with the
system in~\eqref{ODESY} and evaluated at~\eqref{ODE:EQUIL}
presents: 
\begin{itemize}
\item either three real, negative eigenvalues, 
\item red one real, negative eigenvalue
and two complex conjugate eigenvalues with negative real part.
\end{itemize}
There exist parameters that realize each of the cases above.
\medskip

When instead
\begin{equation}\label{LACOSU3}
\left( \alpha f_0 + \gamma + \eta \right) \left( 2 \alpha f_0 \eta + \gamma \eta \right)< \beta \eta \vartheta f_0,\end{equation}
the equilibrium in~\eqref{ODE:EQUIL} is linearly unstable for the system in~\eqref{ODESY}
and the associated Jacobian matrix presents
one real, negative eigenvalue
and two complex conjugate eigenvalues with positive real part.
\medskip

Finally, when
\begin{equation}\label{LACOSU2}
\left( \alpha f_0 + \gamma + \eta \right) \left( 2 \alpha f_0 \eta + \gamma \eta \right)= \beta \eta \vartheta f_0,\end{equation}
the linear stability of the system is undetermined, since the Jacobian matrix
possesses one real, negative eigenvalue and two complex, purely imaginary
eigenvalues.
\end{theorem}

Condition~\eqref{prim:ca} is interesting since, being satisfied for small values of~$\vartheta$,
it suggests that moderate (rather than excessive) controlled burning
could provide more stable outcomes.
We can further determine the relevant factors in making decisions on controlled burning.
In environments with low levels of natural regeneration, 
less oscillations in the prescribed burning program (obtained by reducing $\vartheta$) are desirable, as the factor on the left-hand side of~\eqref{prim:ca} will be lower.
If the regenerative factor~$\gamma$ is high, however, or the optimal bushfire intensity~$f_0$ is high, then condition~\eqref{prim:ca} is satisfied
and 
higher oscillations in the controlled burning program remain permissible without sacrificing stability.

\begin{figure}[h]
    \centering
    \includegraphics[height=4.79cm]{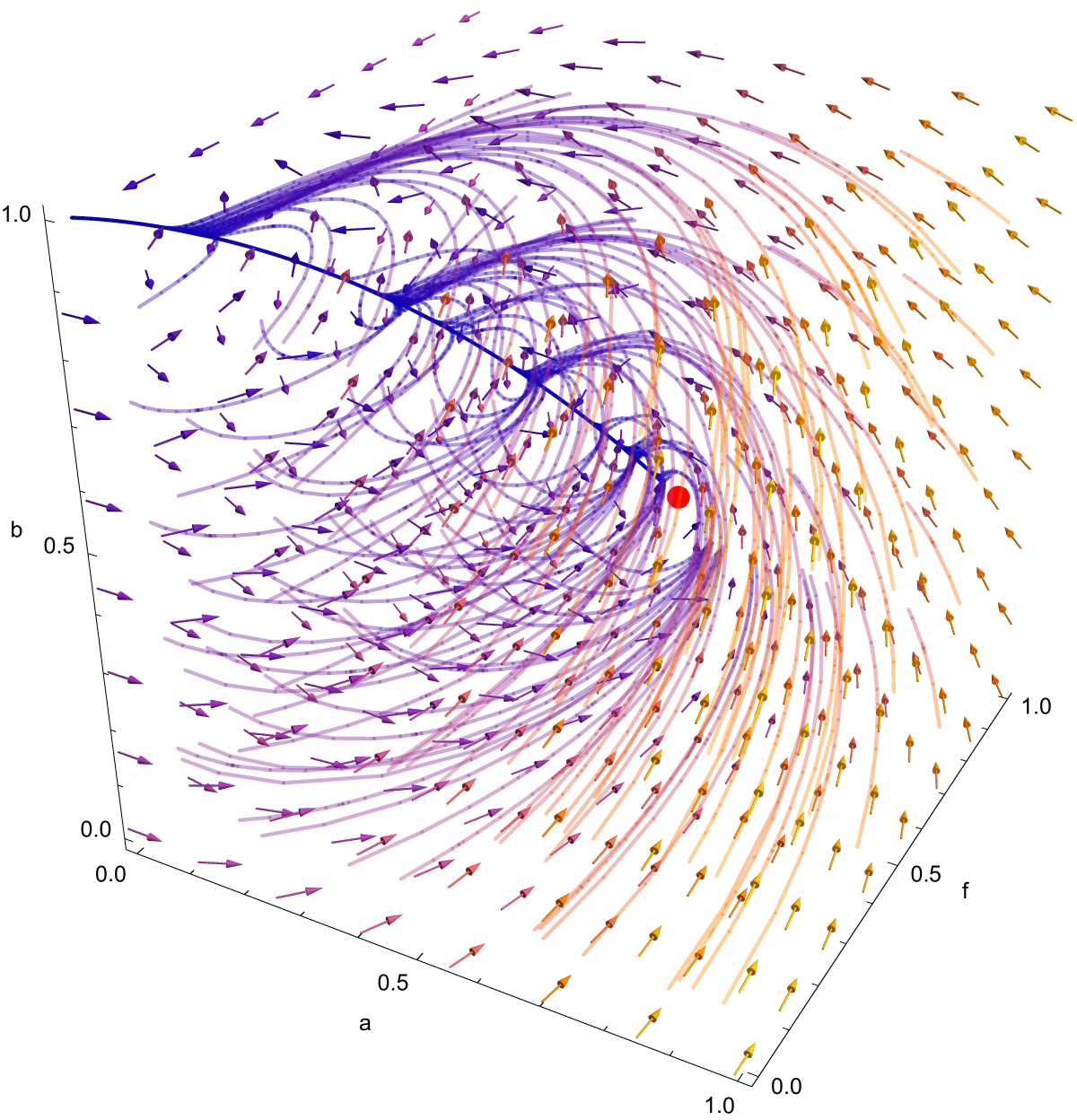}$\quad$\includegraphics[height=4.79cm]{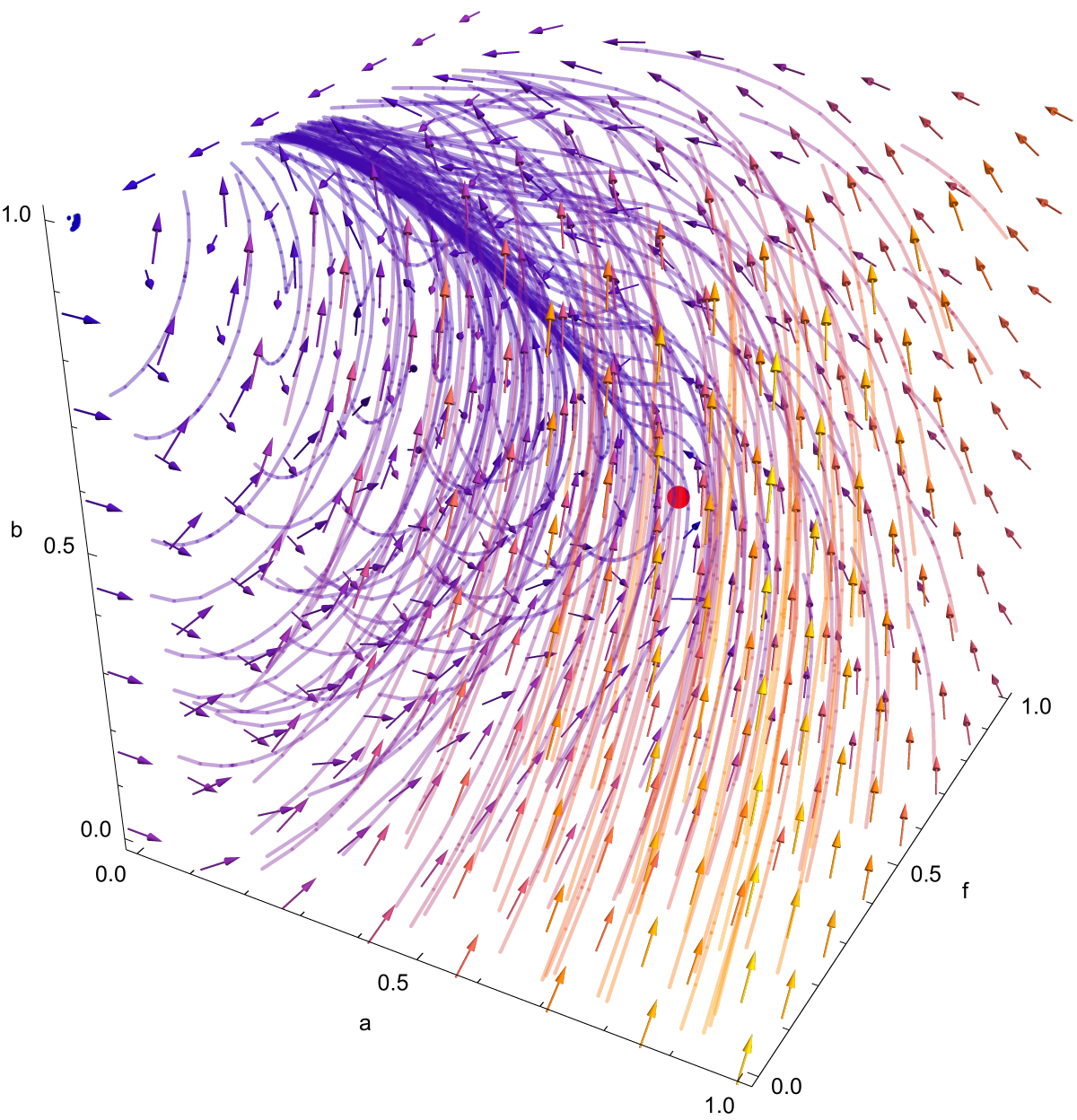}$\quad$\includegraphics[height=4.79cm]{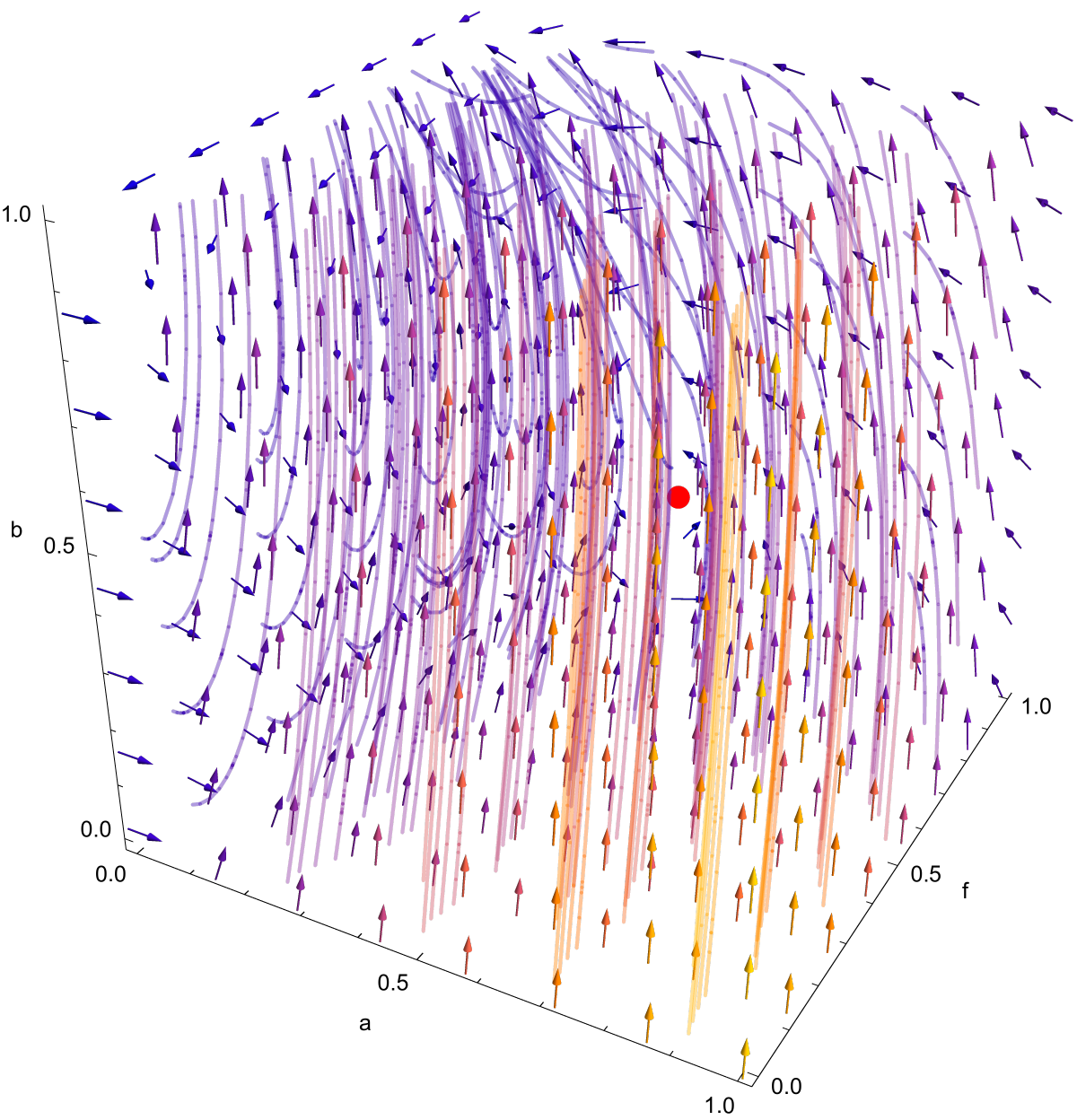}
    \caption{\sl Sketch of the streamlines of~\eqref{ODESY}.
    We adopt here the choice of parameters~$\alpha:=\beta:=\gamma:=\eta:=f_0:=1$ and~$\zeta:=2.5$, which satisfies~\eqref{LACO1} and~\eqref{LACO2}, and~$\vartheta\in\{0.1,\,3,\,20\}$. These values of~$\vartheta$ exemplify, respectively, \eqref{prim:ca}, \eqref{LACOSU2}, and~\eqref{LACOSU3}. The images confirm the transition towards instability.}
    \label{fig1023fgm}
\end{figure}

\begin{figure}[h]
    \centering
    \includegraphics[height=2.79cm]{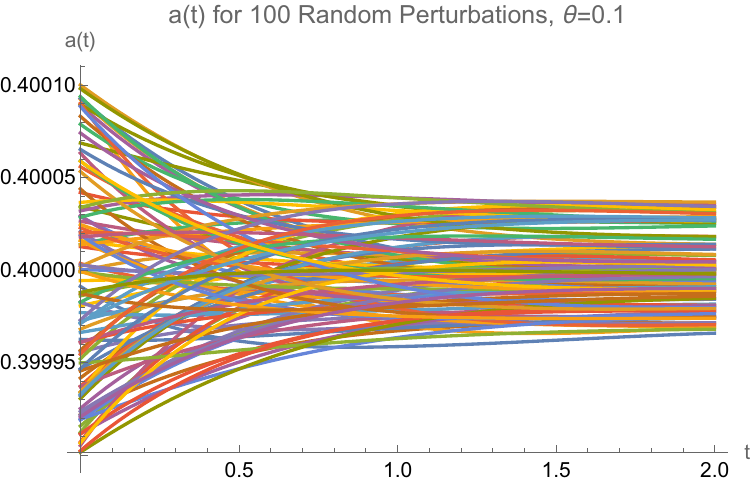}$\quad$\includegraphics[height=2.79cm]{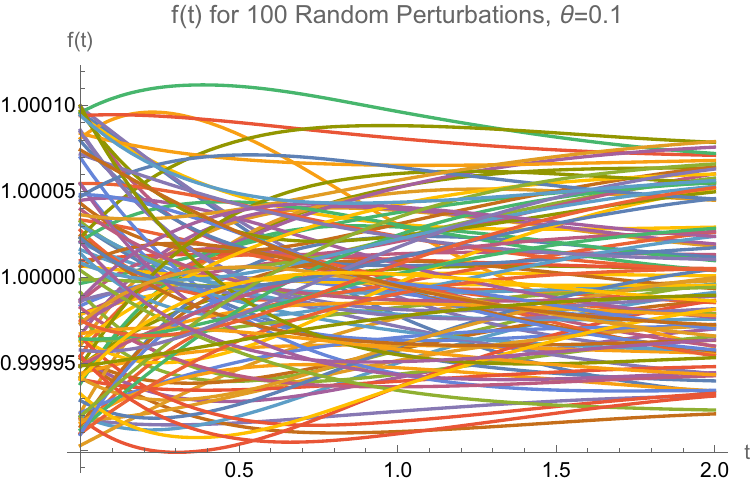}$\quad$\includegraphics[height=2.79cm]{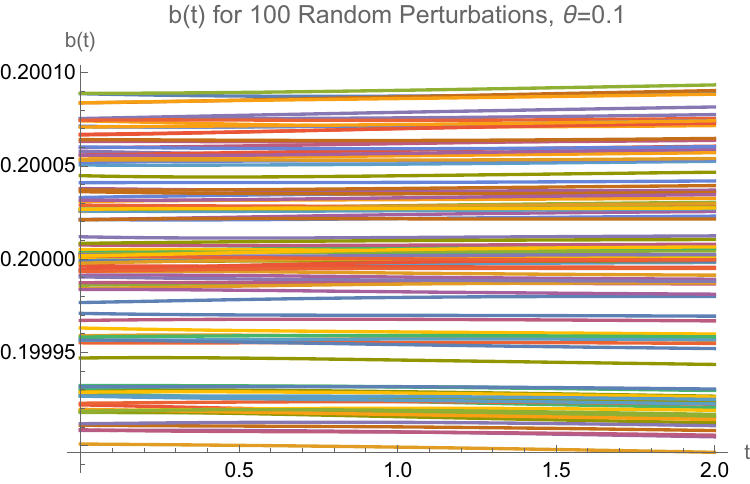}\\
    \includegraphics[height=2.79cm]{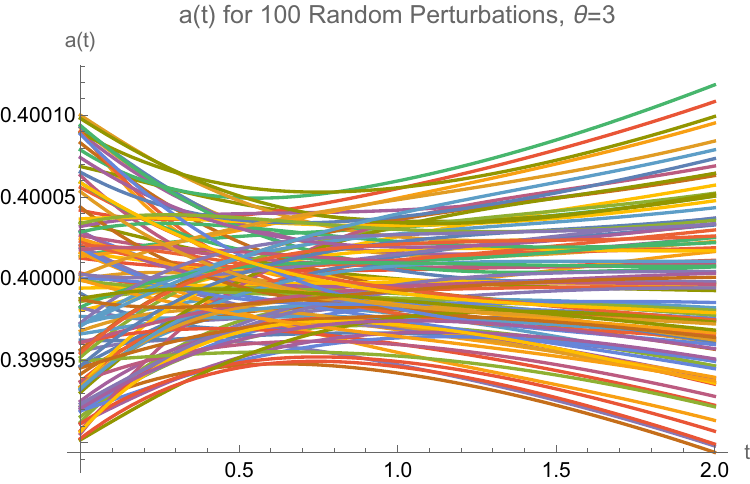}$\quad$\includegraphics[height=2.79cm]{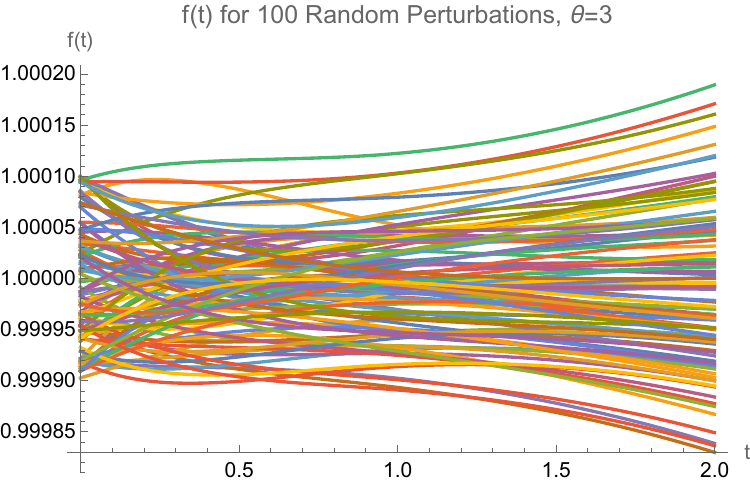}$\quad$\includegraphics[height=2.79cm]{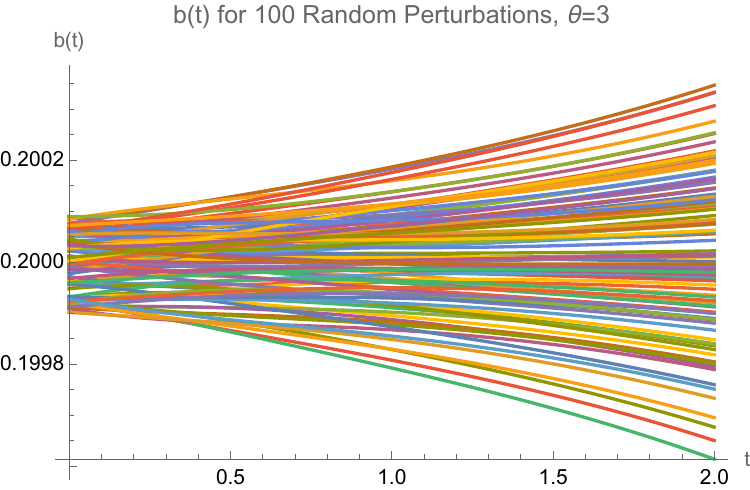}\\
    \includegraphics[height=2.79cm]{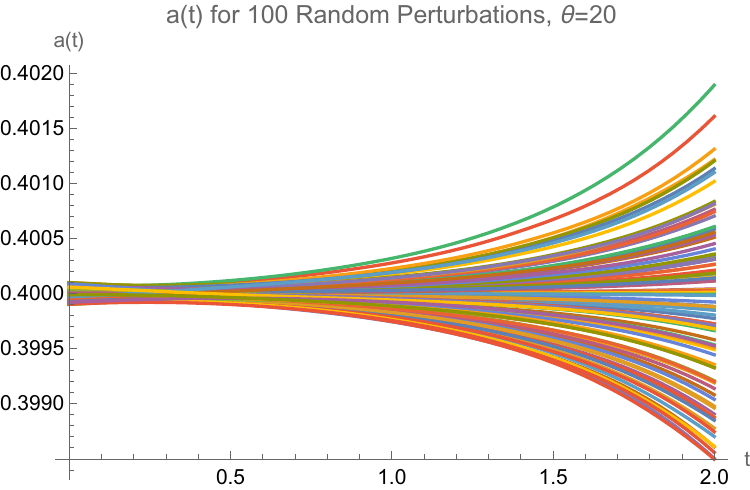}$\quad$\includegraphics[height=2.79cm]{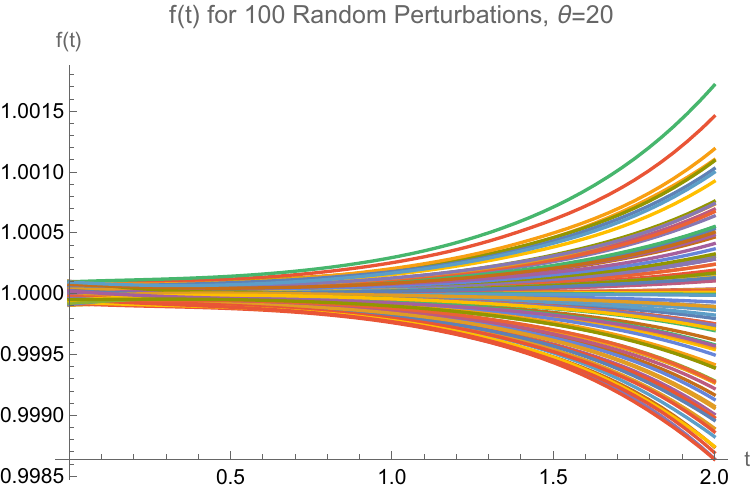}$\quad$\includegraphics[height=2.79cm]{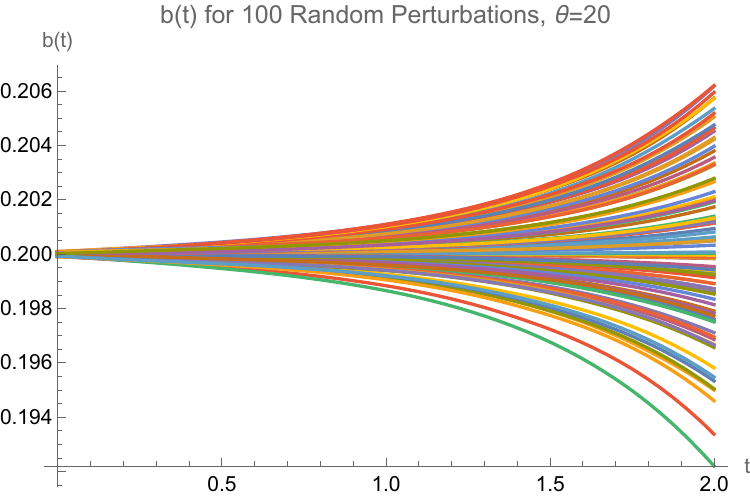}
    \caption{\sl Trajectories of the streamlines of~\eqref{ODESY}
   with the choice of parameters~$\alpha:=\beta:=\gamma:=\eta:=f_0:=1$,
   $\zeta:=2.5$, and~$\vartheta\in\{0.1,\,3,\,20\}$. These values of~$\vartheta$ exemplify, respectively, \eqref{prim:ca}, \eqref{LACOSU2}, and~\eqref{LACOSU3}. The initial points
   are 100 random perturbations of the equilibrium in~\eqref{ODE:EQUIL}.
   The maximal size of the perturbation is~$\pm0.0001$.
   The images confirm the transition towards instability.}
    \label{fig1023fgm.bis}
\end{figure}

The transition to instability, also passing through the critical case in~\eqref{LACOSU2},
is visualized in Figures~\ref{fig1023fgm} and~\ref{fig1023fgm.bis}.

Regarding Theorems~\ref{THSTA} and~\ref{asdfdvfb20rofjgQLNWedfbVSTA},
we observe that, while structural parameters are often difficult to measure in
practice, a key advantage of our analysis is that it relies primarily on a
single parameter, namely~$\theta$, and, crucially, only on its sign. 
In this spirit, although distinguishing between ``reactive'' and ``proactive''
strategies is often impractical---being implicitly tied to the time scale used
to differentiate ``actions'' from ``reactions''---identifying the underlying
blueprint of these ideal programs and assessing their viability through
rigorous, quantitative analysis can provide valuable insights into the
fundamental structure of a complex phenomenon (see Figures~\ref{fig1},
\ref{fig3}, and~\ref{fig2} for several comparisons between natural bushfires
and controlled burning activities which possibly hint at different management
strategies).

\subsection{Further observations}
Regarding the condition~$b\le a$ describing the fact that the land devoted to
preventive burning is less than or equal to the total available land, a
natural question is whether all
solutions of~\eqref{ODESY} fulfill this prescription, e.g. whether or not
when the initial condition satisfies~$b(0)<a(0)$ one has that~$b(t)\le a(t)$,
or at least whether or not this is true under suitable assumptions on the structural
parameters~$\alpha$, $\beta$, $\gamma$, $\zeta$, $\eta$, $\theta$, and~$f_0$. The next result states that this is not the case:

\begin{theorem}\label{0jwdoflvnre:2}
Given any~$\alpha$, $\beta$, $\gamma$, $\zeta$, $\eta$, $\theta$, and~$f_0$ as
above, there exist initial conditions~$a(0)\in(0,1)$,
$b(0)\in\left(0,a(0)\right)$, and~$f(0)\in(0,+\infty)$ and a time~$T>0$ such
that the solution of~\eqref{ODESY} satisfies
 \EQUAZIONE{0< a(T)<b(T)<1.}
\end{theorem}

The situation put forth in Theorem~\ref{0jwdoflvnre:2} is interesting for
practical applications, since it highlights the risk of exhausting the
available land due to (a combined effect of bushfires and) preventive burning.
\medskip

While, according to Theorem~\ref{0jwdoflvnre:2}, the condition~$b<a$ is not
preserved under the flow of~\eqref{ODESY}, the conditions~$f>0$
and~$a\in(0,1)$ do persist (as long as~$0\le b\le a$):

\begin{theorem}\label{0jwdoflvnre:3}
Let~$\alpha$, $\beta$, $\gamma$, $\zeta$, $\eta$, $\theta$, and~$f_0$ be as
above.
Suppose that, for some~$T>0$, the solution of~\eqref{ODESY}
in the interval~$[0,T]$ is such that~$a(t)\in(0,1)$, $f(t)\in(0,+\infty)$, and~$b(t)\in\big[0, a(t)]$
for all~$t\in(0,T)$.

Then, $a(T)\in(0,1)$ and~$f(T)\in(0,+\infty)$.
\end{theorem}

We stress that Theorems~\ref{0jwdoflvnre:2} and~\ref{0jwdoflvnre:3} are valid
both when~$\theta\in(0,+\infty)$ and when~$\theta\in(-\infty,0)$.

We also emphasize that the results in Theorems~\ref{0jwdoflvnre:2} and~\ref{0jwdoflvnre:3} do not imply that the model presented here is unrealistic. On the contrary, Theorem~\ref{0jwdoflvnre:2} should be interpreted as a caution against burning too much land (potentially all the available area)
through control burnings. Meanwhile, Theorem~\ref{0jwdoflvnre:3} serves as a consistency result, provided that the available region is not exhausted.

\subsection{Caveats and perspectives}
We stress that, as is usual when dealing with the potential implications of
theoretical works on concrete scenarios, scientific results must be correctly
interpreted, and a pinch of salt is necessary when attempting to translate
them into actual policies.
On the one hand, one of the advantages of the model~\eqref{ODESY} lies in its
intrinsic simplicity and the relatively small number of structural parameters
involved in describing complex interactions.
This allows for a thorough analytical treatment with explicit quantifications.
On the other hand, the model~\eqref{ODESY} is a simplification of the highly
complex, high-dimensional dynamics involved in bushfire management.
For instance, the equation modeling bushfire frequency in~\eqref{FREQ:EQMD}
can be enriched by incorporating a seasonal term (this can be achieved by
adding a periodic forcing term to its right-hand side).
We wish to also mention two other specific factors.
Latent flammability (see~\cite{LINZYL}) and biodiversity (see~\cite{BIOD}), can each be impacted by hazard reduction, leading to another feedback effect that could influence long-time dynamics.

Moreover, other key factors could be considered in the system~\eqref{ODESY},
such as public awareness, public perception of risks, involvement of local
stakeholders, political trends, media coverage, climate considerations, impact
on the quality of air, alternative managing procedures such as selective
logging and mulching practices, and restoration activities (which themselves
evolve over time and both influence and are influenced by bushfire events).
To study the combination of these factors, it will also be important to
distinguish the different time scales at which they take place (for example,
one can model both regular and catastrophic wildfire events as either random
or periodic, but their frequency differs, and climate change is causing the
period to shrink and the frequency to increase). 
Although incorporating all these factors would likely hinder a fully
analytical treatment of the problem, we plan to revisit these considerations
in a future work.

\subsection{Organisation of the paper}
The rest of this paper is organized as follows.
We give a sketch of the proofs of the main theorems in Section~\ref{SKEE}.
The detailed proofs of these results are presented in
Section~\ref{sqdfvfbqwfegrhtyjetyrsdDSFGHJD}, with the demonstration of some
auxiliary algebraic results deferred to the appendix.
Some conclusions will be summarized at the end of the paper.

\section{Insights into the main proofs}\label{SKEE}
One of the main focuses of this paper is the linear stability and instability analysis of the new model for prescribed burning that we have introduced. The methodology employed is classical, relying on the identification of equilibria and the spectral analysis of the linearized system. In particular, complex eigenvalues of the Jacobian matrix with positive real parts correspond to unstable modes: oscillatory patterns arise when the imaginary part is nonzero, while real, positive eigenvalues indicate instabilities along escaping directions. Conversely, eigenvalues with negative real parts correspond to stable configurations.

The spectral analysis in our case is nontrivial, as the characteristic polynomial of the Jacobian matrix is of third degree. While, in principle, such a polynomial can be explicitly factorized, the resulting expression would be cumbersome and practically intractable, especially due to the presence of several structural parameters. Therefore, to determine the stability or instability of the system, we develop and employ a number of bespoke algebraic arguments.

Specific arguments from the theory of ordinary differential equations will also be used to analyze the trajectories in particular cases of interest.

The mathematical methods employed are elementary in nature, requiring no advanced background. However, they are carefully tailored to suit the model under consideration and to highlight the connection between the theoretical findings and their practical implications for controlled burning planning.

\section{Proofs of the main results}\label{sqdfvfbqwfegrhtyjetyrsdDSFGHJD}

\subsection{Some preliminary algebraic observations} 
The spectral analysis of the linearized system is not standard,
since explicit solutions, though available, would be practically difficult or impossible to manage.
To circumvent this difficulty, we introduce here some auxiliary algebraic results.
Indeed,
in the proof of Theorems~\ref{THSTA} and~\ref{asdfdvfb20rofjgQLNWedfbVSTA},
when establishing that the system may possess, depending on the paramenters,
real eigenvalues as well as complex ones, the following observations will come in handy:

\begin{lemma}\label{ASDEFGRbqrfewg2450t34ughbnasftg34h5ytusdw.0}
Let~$B$, $C$, $D\in(0,+\infty)$, with~$BC\le3 D$.
Then, the polynomial~$\lambda^3+B\lambda^2+C\lambda+D$ possesses two
complex conjugate roots.

If also~$BC<D$, then the real part of these roots is positive.
\end{lemma}

\begin{lemma}\label{ASDEFGRbqrfewg2450t34ughbnasftg34h5ytusdw}
Let \begin{equation}\label{ASDEFGRbqrfewg2450t34ughbnasftg34h5ytusdw2}
B:= \alpha f_0 + \gamma + \eta,\qquad C:= 2 \alpha f_0 \eta + \gamma \eta \qquad{\mbox{and}}\qquad E:=-\beta\eta\theta f_0.\end{equation}
Let also
\begin{equation}\label{qowhfnGBnsdfmNRyasdhfjTsdRV}\triangle:=
B^2C^2-4C^3-4B^3E-27E^2+18BCE.
\end{equation}
Given~$\theta\in\R\setminus\{0\}$, there are suitable choices of the structural parameters~$\alpha$, $\beta$,
$\gamma$, $\zeta$, $\eta$, and~$f_0$, fulfilling
conditions~\eqref{LACO1},
\eqref{LACO2} and $BC > E$, which reads
\begin{equation}\label{CJE7q829erwfgobk-DWFVD}
\left( \alpha f_0 + \gamma + \eta \right) \left( 2 \alpha f_0 \eta + \gamma \eta \right)> -\beta \eta \theta f_0,
\end{equation}
allowing for both of the possibilities~$\triangle>0$
and~$\triangle<0$.\end{lemma}

For the convenience of the reader, the proof of these results is deferred to the appendix.
Notice that the parameter~$\theta$ does not play a role in Lemma~\ref{ASDEFGRbqrfewg2450t34ughbnasftg34h5ytusdw},
allowing us to use the above result both when~$\theta\in(0,+\infty)$
and when~$\theta\in(-\infty,0)$.

In our framework the interest of
Lemma~\ref{ASDEFGRbqrfewg2450t34ughbnasftg34h5ytusdw} is that it detects the
sign of the discriminant~$\triangle$ of the cubic equation
\EQUAZIONE{ \lambda^3 +B\lambda^2 + C\lambda+E=0,}
where~$B$, $C$ and~$E$ are as in~\eqref{ASDEFGRbqrfewg2450t34ughbnasftg34h5ytusdw2}.

Specifically, if the discriminant is positive, the cubic equation has three
distinct real roots, while if the discriminant is negative, the cubic equation
has one real root and two complex conjugate roots (see e.g. Section~6.2
in~\cite{MR1409810} for the role of discriminants in polynomial equations).
Therefore, Lemma~\ref{ASDEFGRbqrfewg2450t34ughbnasftg34h5ytusdw} will allow us
to conclude that both the case of real and complex eigenvalues occur in our
system, according to the statements of Theorems~\ref{THSTA} and~\ref{asdfdvfb20rofjgQLNWedfbVSTA}.

\subsection{Proof of Theorem~\ref{THSTA}}\label{DECASCV:WEFER}
As is customary, the stability of an equilibrium is determined by the Jacobian matrix (at the equilibrium).
In this setting, one says that the equilibrium is \emph{linearly stable} if all
eigenvalues have negative real parts (this, in fact, ensures that the fixed
point is locally attractive with exponential strength, see e.g.  Proposition~6
on page~378 of~\cite{MR698947}) and is \emph{linearly unstable} if at least
one of the eigenvalues has a positive real part.

The Jacobian matrix associated with the system of ordinary differential
equations in~\eqref{ODESY} takes the form
\begin{equation}\label{JACOMAT5} J (a,f,b):=
\begin{bmatrix}
-\alpha f - \gamma & -\alpha a & -\beta \\
\zeta & -\eta & 0 \\
\theta (f-f_0) & \theta a & 0
\end{bmatrix}.\end{equation}

This and~\eqref{ODE:EQUIL} give that
\begin{equation}\label{JACOMAT5CONLOSTA} J(a_\star, f_\star, b_\star) =
\begin{bmatrix}
-\alpha f_0 - \gamma & -\frac{\alpha\eta f_0}{\zeta} & -\beta \\
\zeta & -\eta & 0 \\
0 &  \frac{\eta \theta f_0}{\zeta} & 0
\end{bmatrix}.\end{equation}

The eigenvalues of this matrix are determined by the third-degree
characteristic equation
\begin{equation}\label{FORfo6}
P(\lambda):=
\lambda^3 + \left( \alpha f_0 + \gamma + \eta \right)\lambda^2 + \left( 2 \alpha f_0 \eta + \gamma \eta \right)\lambda - \beta \eta \theta f_0= 0.
\end{equation}
Key for our analysis is that the characteristic polynomial has the form
\begin{equation} \label{FORfo5}\begin{split}&P(\lambda)=\lambda^3 +B\lambda^2 + C\lambda - D,
\\&{\mbox{with~$B$, $D\in(0,+\infty)$ and~$C\in\R$.}}
\end{split}
\end{equation}
In this case we know that~$C\in(0,+\infty)$, but we are not using this
information in what follows.

Since~$P(0)= -D<0$ and~$P(+\infty)=+\infty$,
there exists at least one root~$\lambda_1$ of~$P$ lying in the interval~$(0,+\infty)$,
which proves that the equilibrium in~\eqref{ODE:EQUIL} is linearly unstable.

Let us now consider the other (possibly complex) roots of~$P$, denoted by~$\lambda_2$ and~$\lambda_3$.
Since
\EQUAZIONE{\begin{split}&
\lambda^3 +B\lambda^2 + C\lambda - D=P(\lambda)=(\lambda-\lambda_1)(\lambda-\lambda_2)(\lambda-\lambda_3)\\&\qquad=
\lambda^3 - (\lambda_1+\lambda_2+\lambda_3) \lambda^2 + 
(\lambda_1 \lambda_2+ \lambda_1 \lambda_3+\lambda_2 \lambda_3 )  \lambda  - \lambda_1 \lambda_2 \lambda_3,
\end{split}}
it follows that
\EQUAZIONE{ -B=\lambda_1+\lambda_2+\lambda_3\qquad{\mbox{and}}\qquad
D=\lambda_1 \lambda_2 \lambda_3.}
Since~$\lambda_1\in(0,+\infty)$, we deduce that
\begin{equation}\label{COL:i}\lambda_2 \lambda_3=\frac{D}{\lambda_1}\in(0,+\infty)\end{equation}
and that
\begin{equation}\label{COL:ii} \lambda_2+\lambda_3=-B-\lambda_1\in(-\infty,0).\end{equation}

Now we distinguish two cases: either~$\lambda_2$, $\lambda_3$ are real, or are complex conjugate.
If they are real, it follows from~\eqref{COL:ii} that at least one of them is negative, say~$\lambda_2\in(-\infty,0)$. Hence, by~\eqref{COL:i}, also~$\lambda_3=\frac{D}{\lambda_1\lambda_2}\in(-\infty,0)$.

If instead~$\lambda_2$, $\lambda_3$ are complex conjugate, 
we write~$\lambda_2=x+iy$ and~$\lambda_3=x-iy$, for some~$x\in\R$ and~$y\in\R\setminus\{0\}$,
and we deduce from~\eqref{COL:ii} that~$x\in(-\infty,0)$.

This completes the eigenvalue analysis stated in Theorem~\ref{THSTA}.

Lemma~\ref{ASDEFGRbqrfewg2450t34ughbnasftg34h5ytusdw} guarantees that both of
the above scenarios are possible, depending on the choice of the parameters.
This finishes the proof of Theorem~\ref{THSTA}.~\hfill$\qed$

\subsection{Proof of Theorem~\ref{STAINDO}}
We denote by~$y$ the (column) vector with entries~$a$, $f$, $b$ (the
corresponding equilibrium in~\eqref{ODE:EQUIL} will be denoted by~$y_\star$).
Hence, we rephrase~\eqref{ODESY} in the form
\begin{equation}\label{E3rVOLI0-1}
\dot y=Y(y),\end{equation}
for a suitable (nonlinear) map~$Y:\R^3\to\R^3$. 
In this setting, the equilibrium property of~$y_\star$ reads~$Y(y_\star)=0$.

We use the setting of~\eqref{JACOMAT5CONLOSTA} with the short
notation~$J_\star:=J(a_\star,f_\star,b_\star)$.
By virtue of Theorem~\ref{THSTA}, we consider the eigenvalues~$\lambda_1$,
$\lambda_2$, and~$\lambda_3$ of~$J_\star$, with~$\lambda_1\in(0,+\infty)$
and~$\Re\lambda_2$, $\Re\lambda_3\in(-\infty,0)$.

We perform a Schur triangulation (see e.g. Section~2.3 in~\cite{MR2978290}),
namely we take a unitary matrix $Q$ such that
\EQUAZIONE{ J_\star=Q^{-1}UQ,}
where~$U$ is an upper triangular matrix of the form
\EQUAZIONE{ U=\left(\begin{matrix}
\lambda_1 & * & *\\
0&\lambda_2&*\\
0&0&\lambda_3
\end{matrix}\right).}
Here, ``$*$'' denotes some complex entry (possibly different from one
another).

We consider the change of variable~$x:=Qy$ and define
\EQUAZIONE{ x_\star:=Qy_\star\qquad{\mbox{and}}\qquad
X(x):=QY(Q^{-1}x).}
In this way, equation~\eqref{E3rVOLI0-1} becomes
\begin{equation}\label{E3rVOLI0-2}
\dot x=Q\dot y=QY(y)=QY(Q^{-1}x)=X(x),\end{equation}
and note that~$X(x_\star)=QY(y_\star)=0$.

When needed, we will use coordinates by writing~$x_j$ for the~$j$th entry of
the vector~$x$ and~$X_j$ for the~$j$th entry of the vector~$X$.

In this scenario, we have that~$x_\star$ is an equilibrium for the system of
ordinary differential equations in~\eqref{E3rVOLI0-2}, which is in turn
equivalent to the original system in~\eqref{ODESY} up to a linear mapping.
The corresponding Jacobian matrix satisfies
\begin{equation}\label{0ojfnvt9854nb050325Xe356mAS7} \begin{split}&\frac{\partial X}{\partial x}(x_\star)=
Q\frac{\partial Y}{\partial y}(Q^{-1}x_\star)\,Q^{-1}=
Q\frac{\partial Y}{\partial y}(y_\star)\,Q^{-1}\\&\qquad=QJ_\star\, Q^{-1}=U=\left(\begin{matrix}
\lambda_1 & * & *\\
0&\lambda_2&*\\
0&0&\lambda_3
\end{matrix}\right).\end{split}\end{equation}

To stabilize the equilibrium~$x_\star$ (or equivalently, up to a linear
mapping, to stabilize the original equilibrium~$y_\star$), we propose the
following system of ordinary differential equations which includes now a
feedback function~$\omega:\R\to\R$:
\begin{equation}\label{OSNDYT2c43v7U3tyDIJFJom13fe} 
\begin{cases}\dot x_1=X_1\big(x_1-\tau(\omega-x_{\star,1}),x_2,x_3\big),\\
\dot x_2=X_2(x),\\
\dot x_3=X_3(x),\\
\dot\omega=\sigma(x_1-\omega)
.\end{cases}\end{equation}
Above $\sigma$, $\tau\ge0$ represent feedback parameters.
When~$\tau=0$, the system \eqref{OSNDYT2c43v7U3tyDIJFJom13fe} is equivalent
to the original one.
When~$\sigma$, $\tau>0$, the feedback function~$\omega$ evolves, accounting
for the unstable direction~$x_1$ and acts on its evolution (the stable
directions~$x_2$ and~$x_3$ are not affected by the feedback).

The practical implementation of feedback mechanisms in controlled burns is inherently complex, and necessarily so, given the structural instability of the underlying system. At the core of this approach lies a quantitative analysis of a derived variable, denoted~$x_1$, which captures the aggregate effects of the original variables~$(a, f, b)$, representing available fuel, bushfire intensity, and the extent of controlled burns, respectively.
We remark that this quantity~$x_1$ is not an abstract construct: in fact, once the original variables are measured, $x_1$ is obtained via a simple linear combination. The coefficients involved in this combination are explicitly determined through the above Schur triangulation process, making the formulation of the feedback mechanism fully transparent.

The primary objective of the feedback is to dampen the destabilizing influence of~$x_1$, which is tracked through the evolution of the variable~$\omega$ which quantifies the required adjustment to~$x_1$ in order to restore stability to the system. Translating this adjustment back into the original variables provides an explicit prescription for modifying~$(a, f, b)$ to stabilize bushfire dynamics.

However, this process also highlights a fundamental modeling challenge: in real-world bushfire management, policy interventions can directly influence controlled burns~$b$ and, to some extent, the amount of available fuel
(thus affecting the values of~$a$) through land management practices. In contrast, bushfire intensity~$f$
is largely governed by unpredictable natural factors and thus lies outside direct control.
Despite this limitation, the feedback mechanism described offers a valuable starting point, namely
an explicit tool for guiding interventions. While it may not fully account for all practical constraints, especially regarding the uncontrollable variable~$f$, it serves as a first-order approximation that can be enhanced with heuristic and empirical strategies. Such an approach can be particularly effective when the system operates near an optimal fire activity level~$f_0$, which possibly does not requires major adjustments in the variable~$f$.
Additionally, the explicit expression of the feedback control in certain parameter ranges may suggest indirect actions for policymakers to intervene in bushfire intensity, such as increasing wildfire surveillance, enlarging fire brigades, enhancing early warning systems, improving forest management practices, or allocating more resources to community preparedness and rapid response strategies.

We stress that~$x_{\star\star}:=(x_\star, x_{\star,1})$ is an equilibrium
for~\eqref{OSNDYT2c43v7U3tyDIJFJom13fe} and it follows
from~\eqref{0ojfnvt9854nb050325Xe356mAS7} that the corresponding Jacobian
matrix has the form
\EQUAZIONE{ J_{\star\star}:=\left(\begin{matrix}
\lambda_1 & * & *&-\tau\lambda_1\\
0&\lambda_2&*&0\\
0&0&\lambda_3&0\\
\sigma&0&0&-\sigma
\end{matrix}\right).}
The corresponding eigenvalues are therefore
\EQUAZIONE{\lambda_2, \quad\lambda_3,\quad
\frac{\lambda_1-\sigma  +\sqrt{ (\sigma+\lambda_1)^2  - 4 \sigma \tau \lambda_1 } }2,\quad{\mbox{and}}\quad
\frac{\lambda_1-\sigma  -\sqrt{ (\sigma+\lambda_1)^2  - 4 \sigma \tau \lambda_1 } }2.}
The real part of all these eigenvalues is negative as long
as~$\sigma>\lambda_1$
and~$\tau>\frac{(\sigma+\lambda_1)^2}{4\sigma\lambda_1}$, thus providing the
desired stability property. The proof of Theorem~\ref{STAINDO} is thereby
complete.~\hfill$\qed$

\subsection{Proof of Theorem~\ref{asdfdvfb20rofjgQLNWedfbVSTA}}
The matrix analysis in Section~\ref{DECASCV:WEFER} can be used here
with~$\theta=-\vartheta\in(-\infty,0)$ and accordingly, by~\eqref{FORfo6}, the
characteristic polynomial takes the form
\EQUAZIONE{ P(\lambda):=
\lambda^3 + \left( \alpha f_0 + \gamma + \eta \right)\lambda^2 + \left( 2
\alpha f_0 \eta + \gamma \eta \right)\lambda + \beta \eta \vartheta f_0.}
In particular,
\EQUAZIONE{\begin{split}&P(\lambda)=\lambda^3 +B\lambda^2 + C\lambda + D,
\\&{\mbox{with~$B$, $C$, $D\in(0,+\infty)$.}}
\end{split}}
This is a cubic polynomial and thus it possesses at least one real
root~$\lambda_1$.
Since all the coefficients of~$P$ are positive,
\begin{equation}\label{TUSFVNEFVGA2R}
{\mbox{all real roots must be negative.}}\end{equation}

Now, we distinguish three cases: either
\begin{equation}\label{CASO-diwsdnt1}
\left( \alpha f_0 + \gamma + \eta \right) \left( 2 \alpha f_0 \eta + \gamma \eta \right)> \beta \eta \vartheta f_0,\end{equation}
or
\begin{equation}\label{CASO-diwsdnt2}
\left( \alpha f_0 + \gamma + \eta \right) \left( 2 \alpha f_0 \eta + \gamma \eta \right)< \beta \eta \vartheta f_0,\end{equation}
or
\begin{equation}\label{CASO-diwsdnt3}
\left( \alpha f_0 + \gamma + \eta \right) \left( 2 \alpha f_0 \eta + \gamma \eta \right)= \beta \eta \vartheta f_0.\end{equation}

We show that
\begin{equation}\label{Cw4fr6523ed27902we3wrtO-diwsdnt22er3tw4rgeh}
{\mbox{when~\eqref{CASO-diwsdnt1} holds true, 
the equilibrium is linearly stable.}}\end{equation}
To this end, in view of~\eqref{TUSFVNEFVGA2R}, it suffices to consider the
case in which~$P$ presents two complex conjugate roots~$\lambda_2=x+iy$
and~$\lambda_3=x-iy$, with~$x\in\R$ and~$y\in\R\setminus\{0\}$.

We can rewrite~\eqref{CASO-diwsdnt1} as~$BC>D$ and therefore
\EQUAZIONE{ P(-B)=-B^3 +B^3 -B C + D<0.}
Since~$P(+\infty)=+\infty$, this gives that
\begin{equation}\label{CA-di213e3rwegbqsdafvfdwsdnt1}
\lambda_1\in(-B,0).\end{equation}

Now, since
\begin{equation}\label{asdcjhvd9iSk20owqfjevlb0pwqfv934tgjhHSIDK2}\begin{split}&
\lambda^3 +B\lambda^2 + C\lambda + D=P(\lambda)=(\lambda-\lambda_1)(\lambda-\lambda_2)(\lambda-\lambda_3)\\&\qquad=
\lambda^3 - (\lambda_1+\lambda_2+\lambda_3) \lambda^2 + 
(\lambda_1 \lambda_2+ \lambda_1 \lambda_3+\lambda_2 \lambda_3 )  \lambda  - \lambda_1 \lambda_2 \lambda_3
\end{split}\end{equation}
we thus have
\begin{equation}\label{CA-di213e3rwegbqsdafvfdwsdnt11} 2x=\lambda_2+\lambda_3=\lambda_1+\lambda_2+\lambda_3-\lambda_1=
-B-\lambda_1.\end{equation}
It follows from this and~\eqref{CA-di213e3rwegbqsdafvfdwsdnt1} that~$x<0$,
which completes the proof
of~\eqref{Cw4fr6523ed27902we3wrtO-diwsdnt22er3tw4rgeh}.

Further, by virtue of Lemma~\ref{ASDEFGRbqrfewg2450t34ughbnasftg34h5ytusdw},
we have that both the case of three real roots
and that of one real root and two complex conjugate roots can occur.
This completes the case in which~\eqref{CASO-diwsdnt1}
is satisfied.

Let us now consider the case in which~\eqref{CASO-diwsdnt2} holds true.
In this situation we have that~$BC< D$ and thus, by~\eqref{TUSFVNEFVGA2R}
and Lemma~\ref{ASDEFGRbqrfewg2450t34ughbnasftg34h5ytusdw.0}, we know that
the characteristic polynomial presents a real, negative root
and
two complex conjugate roots with positive real part, also giving that the equilibrium is linearly unstable.

Finally, if~\eqref{CASO-diwsdnt3} holds true, the linear stability of the system is undetermined, since in this case~$BC=D$ and thus
\EQUAZIONE{P(\lambda)=\lambda^3 +B\lambda^2 + C\lambda + BC=(\lambda+B)(\lambda^2+C),}
which presents two purely imaginary roots.
This finishes the proof of Theorem~\ref{asdfdvfb20rofjgQLNWedfbVSTA}.~\hfill$\qed$

\subsection{Proof of Theorem~\ref{0jwdoflvnre:2}}
Consider the initial data~$a(0):=b(0)\in\left(\frac{\gamma}{\alpha
f_0+\beta+\gamma},1\right)$ and~$f(0):=f_0$.
Then, we deduce from~\eqref{ODESY} that
\EQUAZIONE{\begin{split}& \dot b(0)-\dot a(0) =
\theta (f(0)-f_0) a(0)
+\alpha f(0) a(0)+\beta b(0) -\gamma (1-a(0))\\&\quad=\big(
\alpha f(0) +\beta +\gamma\big)a(0) -\gamma >0.\end{split}}
As a result, there exists~$\e>0$ small enough such
that~$b(\e)-a(\e)>0>b(-\e)-a(-\e)$, with~$a(\e)$, $b(\e)\in(0,1)$.
The desired result now follows via a time translation~$t\mapsto t+\e$, by
choosing~$T:=2\e$.~\hfill$\qed$

\subsection{Proof of Theorem~\ref{0jwdoflvnre:3}}
We argue for the sake of contradiction.
Three cases can occur: either~$a(T)=0$, $a(T)=1$, or~$f(T)=0$.

If~$a(T)=0$ for some~$T>0$, then~$b(T)\le a(T)=0$ and thus~$ \dot a(T) =
-\beta b(T) +\gamma >0$.
But then, for~$\e>0$ small enough, $a(T-\e)<a(T)=0$, in contradiction with our
assumptions.

Similarly, if~$f(T)=0$ for some~$T>0$, then~$\dot f(T)=\zeta a(T)>0$ (since we
know already that~$a(T)>0$). But then, for~$\e>0$ small enough,
$f(T-\e)<f(T)=0$, which is a contradiction.

If instead~$a(T)=1$, then~$ \dot a (T)= -\alpha f (T)-\beta b(T) \le -\alpha f
(T)<0$ (since we know already that~$f(T)>0$).
Hence,  for~$\e>0$ small enough, $a(T-\e)>a(T)=1$, which is
absurd.
This finishes the proof.~\hfill$\qed$

\begin{appendix}

\section{}

\subsection{Proof of Lemma~\ref{ASDEFGRbqrfewg2450t34ughbnasftg34h5ytusdw.0}}
Let~$P(\lambda):=\lambda^3+B\lambda^2+C\lambda+D$. 
We claim that
\begin{equation}\label{Sac0jvorhrty}
{\mbox{$P$ possesses two complex conjugate roots.}}
\end{equation}
To check this, we point out that,
since~$P(\lambda)\ge D>0$ when~$\lambda\ge0$,
we have that all the real roots of~$P$ are negative.

Hence, if, for the sake of contradiction, the claim in~\eqref{Sac0jvorhrty}
were false, then~$P$ would have three real roots, say~$-a$, $-b$ and~$-c$, with~$a$, $b$, $c\in(0,+\infty)$.

As a consequence,
\EQUAZIONE{
P(\lambda)=(\lambda+a)(\lambda+b)(\lambda+c)=
\lambda^3+(a+b+c)\lambda^2+(ab+ac+bc)\lambda+abc}
and accordingly~$B=a+b+c$, $C=ab+ac+bc$ and~$D=abc$.

But then~$ BC=(a+b+c)(ab+ac+bc)>3D$, in contradiction with our assumption.

Having established~\eqref{Sac0jvorhrty}, to complete the proof of
Lemma~\ref{ASDEFGRbqrfewg2450t34ughbnasftg34h5ytusdw.0} we need to
show that when~$BC<D$ the real part of these complex roots is positive.

To this end, from the above analysis we know that the three roots of~$P$
have the form~$-a$, $\alpha+i\beta$ and~$\alpha-i\beta$, with~$a\in(0,+\infty)$, $\alpha\in\R$, and~$\beta\in\R\setminus\{0\}$, and our goal is to show that
\begin{equation}\label{alertenv19wefqw-df0gvbhoi}
\alpha\in(0,+\infty).
\end{equation}

For this, we observe that
\EQUAZIONE{\begin{split}
P(\lambda)&=(\lambda+a)(\lambda-\alpha-i\beta)(\lambda-\alpha+i\beta)\\&=
\lambda^3+
(a - 2 \alpha) \lambda^2
+ (\alpha^2 - 2 a \alpha + \beta^2) \lambda
+ a (\alpha^2  + \beta^2 ).
\end{split}}
As a result,
\EQUAZIONE{B=a - 2 \alpha,\quad C=
\alpha^2 - 2 a \alpha + \beta^2\quad{\mbox{and}}\quad
D= a (\alpha^2  + \beta^2 ).}
Consequently,
\EQUAZIONE{\begin{split}
0&<D-BC\\
&=a (\alpha^2  + \beta^2 )-(a - 2 \alpha)(\alpha^2 - 2 a \alpha + \beta^2)\\
&=2\alpha(a^2-2a\alpha+\alpha^2+\beta^2)\\&=
2\alpha\big( (a-\alpha)^2+\beta^2\big),
\end{split}}
from which~\eqref{alertenv19wefqw-df0gvbhoi} plainly follows.~\hfill$\qed$

\subsection{Proof of Lemma~\ref{ASDEFGRbqrfewg2450t34ughbnasftg34h5ytusdw}}
We choose
\begin{equation}\label{FG:QDSCcl2r93tg0:-23erf1} f_0:=1, \qquad\gamma:=\eta:=c\alpha\qquad{\mbox{ and }}\qquad\zeta:=(1+c)\alpha, \end{equation}
with~$c\in(0,+\infty)$.

\begin{figure}[h]
    \centering
    \includegraphics[height=4.79cm]{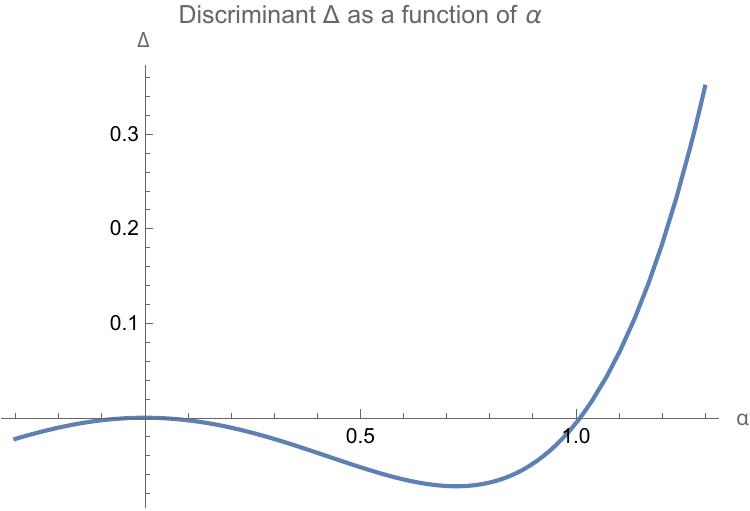}
    \caption{\sl Plot of~$\triangle$ as a function of~$\alpha$ when~$f_0:=1$, $\gamma:=\eta:=c\alpha$, $\zeta:=(1+c)\alpha$, $\beta:=1$, and~$c:=0.1$. Notice that~$\triangle>0$ for large values of~$\alpha$.}
    \label{fig1023fgm.TT1}
\end{figure}

In this way, conditions~\eqref{LACO1} and~\eqref{LACO2} are satisfied and we
see that
\EQUAZIONE{B= (1+2c)\alpha ,\qquad C= c(2  + c)\alpha^2 \qquad{\mbox{and}}\qquad E=-c\alpha\beta\theta,}
with~$\beta\in(0,+\infty)$ and~$\theta\in\R\setminus\{0\}$.

We now construct a case in which~$\triangle>0$.
For this, we suppose that~$\alpha\to+\infty$ (for given~$\beta$, $\theta$,
and~$c$) and we observe that~\eqref{CJE7q829erwfgobk-DWFVD} is fulfilled.
In this situation, we take~$c$ sufficiently small such that
\EQUAZIONE{ \Xi_0:=(1+2c)^2 (2  + c)^2-4c(2  + c)^3>0.}
Then, we infer from~\eqref{qowhfnGBnsdfmNRyasdhfjTsdRV} that
\EQUAZIONE{ \triangle=
(1+2c)^2 c^2(2  + c)^2\alpha^6-4c^3(2  + c)^3\alpha^6
+O(\alpha^4)=
c^2\Xi\alpha^6,}
where
\EQUAZIONE{ \Xi:=
(1+2c)^2 (2  + c)^2-4c(2  + c)^3+O\left(\frac1{\alpha^2}\right)=
\Xi_0+O\left(\frac1{\alpha^2}\right).}
Hence, for~$\alpha$ large enough (possibly depending on~$\beta$, $\theta$,
and~$c$), we conclude that~$\Xi\ge\frac{\Xi_0}2>0$ and
accordingly~$\triangle\ge\frac{c^2\Xi_0\alpha^6}2>0$
(see Figure~\ref{fig1023fgm.TT1} for a specific visualization of this argument).

Now we show that~$\triangle<0$ is also possible. 
{For this, we distinguish two cases, depending on the sign of~$\theta$.

If~$\theta\in(-\infty,0)$,
we consider~$\alpha:=1$,
$\beta:=\frac{1}{|\theta|}$, and~$c\in(0,+\infty)$
sufficiently small.
In this situation, we note that
\EQUAZIONE{\begin{split}
&\left( \alpha f_0 + \gamma + \eta \right) \left( 2 \alpha f_0 \eta + \gamma \eta \right)-\beta \eta |\theta| f_0
= (1+2c) c(2  + c)-c\\&\qquad=c\big((1+2c) (2  + c)-1\big)>c>0
\end{split}}
and therefore~\eqref{CJE7q829erwfgobk-DWFVD} is satisfied.

We also have that~$B= 1+2c$, $C= c(2  + c)$ and~$ E=c$. Thus,
we infer from~\eqref{qowhfnGBnsdfmNRyasdhfjTsdRV} that
\EQUAZIONE{\triangle=-4(1+2c)^3 c+O(c^2)<0,}
as long as~$c$ is sufficiently small
(see Figure~\ref{fig1023fgm.TT2} for a specific visualization of this argument).

\begin{figure}[h]
    \centering
    \includegraphics[height=4.79cm]{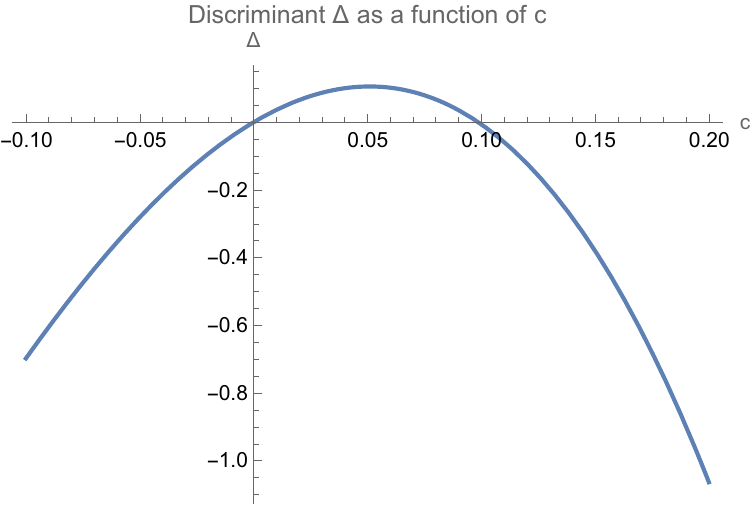}
    \caption{\sl Plot of~$\triangle$ when~$f_0:=1$, $\gamma:=\eta:=c\alpha$, $\zeta:=(1+c)\alpha$, $\alpha:=1$, $\theta:=\beta:=1$.
Notice that~$\triangle>0$ for small values of~$c>0$.}
    \label{fig1023fgm.TT2}
\end{figure}

If instead~$\theta\in(0,+\infty)$, then~\eqref{CJE7q829erwfgobk-DWFVD} is automatically
satisfied. Thus, we consider~$\alpha\searrow0$ (for given~$\beta$,
$\theta$, and~$c$) and we infer from~\eqref{qowhfnGBnsdfmNRyasdhfjTsdRV} that
\EQUAZIONE{\triangle=-27c^2\alpha^2\beta^2\theta^2+O(\alpha^4)<0.}}
These cases complete the proof of
Lemma~\ref{ASDEFGRbqrfewg2450t34ughbnasftg34h5ytusdw}.~\hfill$\qed$
\end{appendix}

\section*{Conclusions}
We have introduced a model comprising a system of three ordinary differential
equations to describe the interplay of bushfires events and
prescribed burning activities in terms of available land.

By simple analytical methods, we have demonstrated that the model possesses a
unique fixed point, representing a steady-state equilibrium of the system.
The stability of this equilibrium is sensitive to the policy governing
prescribed burning.
Specifically:
\begin{itemize}\item
{\em Unstable Equilibrium with Incremental Prescribed Burning (Reactive
Policy):}
When prescribed burning activities are increased in response to bushfire
events, the equilibrium is unstable. 
\item {\em Stable Equilibrium with Prescribed Burning Replacing Natural
Bushfires (Proactive Policy):}
In contrast, when moderate (not excessive) prescribed burning replaces natural bushfires, the equilibrium is stable and the system can maintain a balanced
state. The stability of the equilibrium is also enhanced in the presence of a high regeneration factor.
\end{itemize}

We have also shown that the risks entailed in the Reactive Policy can be
mitigated via a linear feedback control on the system.
\medskip

Though our rigorous analytic results show how consistent policies are
effective in bringing the system into a stable point, these results should not
be seen as the definitive answer to the issue of prescribed burning.
Many specific factors---such as climate conditions, vegetation types, fire
history, and socio-economic implications---must be taken into account in any
forward-thinking policymaking.
However, the stability analysis presented here suggests that policies driven
solely by reactions to bushfire events might be ineffective, potentially
leading to instability. 
Therefore, despite their undeniable efficacy in raising public awareness, to
ensure full effectiveness reactive policies should incorporate feedback
control mechanisms.
On the other hand, proactive policies that focus on planned and controlled
prescribed burning, in moderate quantities and depending on recent bushfire intensity, could be more
effective in prevention and long-term stability.


\begin{bibdiv}
\begin{biblist}

\bib{iwfkhvgri0eotAsd}{article}{
author={Addington, Robert N.},
author={Hudson, Stephen J.},
author={Hiers, J. Kevin},
author={Hurteau, Matthew D.},
author={Hutcherson, Thomas F.},
author={Matusick, George},
author={Parker, James M.}, year={2015},
title={Relationships among wildfire, prescribed fire, and drought in a fire-prone landscape in the south-eastern United States}, journal={Intern. J. Wildland Fire}, number={24}, pages={778--783},
doi={10.1071/WF14187},
url={https://www.publish.csiro.au/wf/wf14187},}

\bib{amo22}{article}{
journal={PLoS One}, year={2022}, volume={17}, issue={5},
doi={10.1371/journal.pone.0240271},
title={Fire use practices, knowledge and perceptions in a West African savanna parkland},
author={Amoako, Esther Ekua},
author={Gambiza, James},
pages={1--25},}

\bib{MR1409810}{book}{
   author={Barbeau, E. J.},
   title={Polynomials},
   series={Problem Books in Mathematics},
   note={Corrected reprint of the 1989 original},
   publisher={Springer-Verlag, New York},
   date={1995},
   pages={xxii+455},
   isbn={0-387-96919-5},
   review={\MR{1409810}},
}

\bib{BEER}{article}{
journal={Math. Computer Modelling},
Volume={13}, Issue={12}, year={1990}, Pages={49--56},
title={The Australian National Bushfire model project},
Author={Beer, Tom},}

\bib{ojdsacnpreoogub987381i2kwertigf}{article}{
title={Hazard burns could be fuelling the fire danger. Indigenous practices offer another way},
year={2023},
author={Butler, Gavin},
url={https://www.sbs.com.au/news/article/hazard-burns-could-be-fuelling-the-fire-danger-indigenous-practices-offer-another-way/pevprmc3f},
journal={SBS News},
}

\bib{BUTZ09}{article}{
author={Butz, Ramona J.},
title={Traditional fire management: historical fire regimes
and land use change in pastoral East Africa},
journal={Intern. J. Wildland Fire}, year={2009}, volume={18}, pages={442--450},}

\bib{CLARKE}{article}{
author={Clarke, Hamish},
author={Cirulis, Brett},
author={Penman, Trent},
author={Price, Owen},
author={Boer, Matthias M.},
author={Bradstock, Ross},
year={2022},
title={The 2019–2020 Australian forest fires are a harbinger of decreased prescribed burning effectiveness under rising extreme conditions},
journal={Scient. Rep.},
VOLUME={12},
ISSUE={1},
URL={https://doi.org/10.1038/s41598-022-15262-y},
pages={1--10},
DOI={10.1038/s41598-022-15262-y},}

\bib{iadjhfw0ireuy03245-63FSp12-35oyu}{article}{
author={Cogos, Sarah},
author={Roturier, Samuel},
author={\"Ostlund, Lars},
year={2020},
title={The origins of prescribed burning in Scandinavian forestry: the seminal role of Joel Wretlind in the management of fire-dependent forests},
JOURNAL={Europ. J. Forest Research},
pages={393--406},
Volume={139},
ISSUE={3},
URL={https://doi.org/10.1007/s10342-019-01247-6},
DOI={10.1007/s10342-019-01247-6},}

\bib{MR4772545}{article}{
   author={Dipierro, Serena},
   author={Valdinoci, Enrico},
   author={Wheeler, Glen},
   author={Wheeler, Valentina-Mira},
   title={A simple but effective bushfire model: analysis and real-time
   simulations},
   journal={SIAM J. Appl. Math.},
   volume={84},
   date={2024},
   number={4},
   pages={1504--1514},
   issn={0036-1399},
   review={\MR{4772545}},
   doi={10.1137/24M1644596},
}

\bib{qa0DOFHWP3REP9H203EPFL9yusdhfiovlf30er-21}{article}{
title={Catastrophic Bushfires, Indigenous Fire Knowledge and Reframing Science in Southeast Australia},
author={Fletcher, Michael-Shawn},
author={Romano, Anthony},
author={Connor, Simon},
author={Mariani, Michela},
author={Yoshi Maezumi, Shira},
journal={Fire}, year={2021}, volume={4}, issue={3}, url={https://doi.org/10.3390/fire4030061},
doi={10.3390/fire4030061}, pages={1--11},}

\bib{MR698947}{book}{
   author={Gallavotti, Giovanni},
   title={The elements of mechanics},
   series={Texts and Monographs in Physics},
   note={Translated from the Italian},
   publisher={Springer-Verlag, New York},
   date={1983},
   pages={xiv+575},
   isbn={0-387-11753-9},
   review={\MR{698947}},
   doi={10.1007/978-3-662-00731-0},
}

\bib{Haikerwal04052015}{article}{
author = {Haikerwal, Anjali}, 
author = {Reisen, Fabienne}, 
author = {Sim, Malcolm R.}, 
author = {Abramson, Michael J.},
author = {Meyer, Carl P.}, 
author = {Johnston, Fay H.}, 
author = {Dennekamp, Martine},
title = {Impact of smoke from prescribed burning: Is it a public health concern?},
journal = {J. Air Waste Manag. Assoc.},
volume = {65},
number = {5},
pages = {592--598},
year = {2015},
doi = {10.1080/10962247.2015.1032445},
URL = {https://doi.org/10.1080/10962247.2015.1032445},
}

\bib{BIOD}{article}{
  title={Could biodiversity loss have increased Australia's bushfire threat?},
  author={Hayward, Matt W},
  author={Ward-Fear, Georgia},
  author={L'Hotellier, Felicity},
  author={Herman, Kerryn},
  author={Kabat, Alexander P},
  author={Gibbons, James P},
  journal={Animal Conservation},
  volume={19},
  number={6},
  pages={490--497},
  year={2016},
  publisher={Wiley Online Library}
}

\bib{MR2144536}{book}{
   author={Hirsch, Morris W.},
   author={Smale, Stephen},
   author={Devaney, Robert L.},
   title={Differential equations, dynamical systems, and an introduction to
   chaos},
   series={Pure and Applied Mathematics (Amsterdam)},
   volume={60},
   edition={2},
   publisher={Elsevier/Academic Press, Amsterdam},
   date={2004},
   pages={xiv+417},
   isbn={0-12-349703-5},
   review={\MR{2144536}},
}

\bib{MR2978290}{book}{
   author={Horn, Roger A.},
   author={Johnson, Charles R.},
   title={Matrix analysis},
   edition={2},
   publisher={Cambridge University Press, Cambridge},
   date={2013},
   pages={xviii+643},
   isbn={978-0-521-54823-6},
   review={\MR{2978290}},
}

\bib{m2edia}{article}{
author={Knight, Ben},
year={2019},
title={Former fire chief calls for more planned burns as fuel loads reach Black Saturday levels},
url={https://www.abc.net.au/news/2019-02-07/black-saturday-fire-fuel-threat-planned-burns-needed/10787050},
journal={www.abc.net.au/news},}

\bib{LAMBE}{collection}{
author={Lambert, Bernard},
title={The French Prescribed Burning Network and
its Professional Team in Pyr\'en\'ees Orientales:
Lessons Drawn from 20 Years of Experience},
series={Best Practices of Fire Use --
Prescribed Burning and Suppression
Fire Programmes in Selected
Case-Study Regions in Europe},
Editor={Montiel, Cristina}, Editor={Kraus, Daniel},
publisher={European Forest Institute},
ISBN={978-952-5453-69-0, 978-952-5453-70-6},
ISSN={1238-8785},
year={2010},
pages={viii+169},
}

\bib{LINZYL}{article}{
  title={Identifying and managing disturbance-stimulated flammability in woody ecosystems},
  author={Lindenmayer, D.},
  author={Zylstra, P.},
  journal={Biological Reviews},
  volume={99},
  number={3},
  pages={699--714},
  year={2024},
  publisher={Wiley Online Library}
}

\bib{2021REP}{article}{
author={Melvin, Mark A.},
title={2021 National
Prescribed Fire Use
Survey
Report},
journal={Technical Report 01-22},
pages={20},
year={2022},
}

\bib{SPARK}{book}{
author={Miller, Claire}, 
author={Hilton, James}, 
author={Sullivan, Andrew}, 
author={Prakash, Mahesh}, 
title={SPARK -- A Bushfire Spread Prediction Tool},
journal={Environmental Software Systems. Infrastructures, Services and Application},
series={Advances in Information and Communication Technology}, publisher={Springer, Cham},
doi={10.1007/978-3-319-15994-2_26},}

\bib{MR1908418}{book}{
   author={Murray, J. D.},
   title={Mathematical biology. I},
   series={Interdisciplinary Applied Mathematics},
   volume={17},
   edition={3},
   note={An introduction},
   publisher={Springer-Verlag, New York},
   date={2002},
   pages={xxiv+551},
   isbn={0-387-95223-3},
   review={\MR{1908418}},
}

\bib{REPO}{book}{
 url={https://knowledge.aidr.org.au/media/4893/overview-of-prescribed-burning-in-australasia.pdf},
 title={National Burning Project.
 Overview of Prescribed Burning
in Australasia.
Report for National Burning Project: Sub-Project 1},
year={2015},
author={Poynter, Mark},
Editor={Kington, Wayne}, 
publisher={Australasian Fire and Emergency Service Authorities Council Limited},
pages={84},
ISBN={978-0-9872065-8-9},}

\bib{NEA1qeA}{article}{
title={Bushfire simulators and analysis in Australia: insights into an emerging sociotechnical practice},
author={Neale, Timothy},
Pages={200--218},
journal={Environm. Hazards},
Volume={17}, year={2018},
Issue={3},}

\bib{q-iw0pf294rf8vj30276jgrth}{article}{
author={Riebold, R. J.}, year={1971}, title={The early history of wildfires and prescribed burning}, journal={Proc. Prescribed Burning Symposium. Charleston, SC},
pages={11--20},}

\bib{RECOMM}{article}{
url={http://royalcommission.vic.gov.au/finaldocuments/summary/PF/VBRC_Summary_PF.pdf},
year={2010},
author={Teague, Bernard},
author={McLeod, Ronald},
author={Pascoe, Susan}, title={2009 Victorian Bushfires Royal Commission. Final Report},
pages={1--42},}

\end{biblist}
\end{bibdiv}
\end{document}